\documentstyle[11pt]{article}

\newcommand{\pall}{ {Part} (2p) }
\newcommand{\palln}{ {Part} (n) }
\newcommand{\pset}{ \{ 1, \ldots , p \} }
\newcommand{\ppset}{ \{ 1, \ldots , 2p \} }
\newcommand{\pp}{ {\cal P}(2p) }
\newcommand{\symmp}{ {\cal S}_p }
\newcommand{\perm}{ {\cal S}_p }
\newcommand{\symmn}{ {\cal S}_n }

\newcommand{\orpi}{ \pi, \varepsilon , \sigma }
\newcommand{\ambda}{ ker }

\newcommand{\ee}{ \varepsilon }
\newcommand{\cs}{ c_{r_1}^{\ee (1)} \cdots  c_{r_n}^{\ee (n)} }
\newcommand{\ls}{ l_{t_{1}}^{\theta_{1}} \cdots  l_{t_{n}}^{\theta_{n}} }
\newcommand{\rs}{ r_1 , \ldots , r_n \in \{ 1, \ldots , s \} }
\newcommand{\ts}{ t_{1}, \ldots , t_{n} \in \{ 1, \ldots , 2s \} }
\newcommand{\tthetas}{ ( t_{1}, \ldots , t_{n} ; \theta_{1}, \ldots ,
\theta_{n} ) }
\newcommand{\rees}{ ( r_1 , \ldots , r_n ; \ee (1), \ldots , \ee (n) ) }
\newcommand{\ees}{ \ee (1), \ldots , \ee (n) \in \{ 1, * \} }
\newcommand{\thetas}{ \theta_{1}, \ldots , \theta_{n} \in \{ 1, * \} }
\newcommand{\pree}{ {\cal P} ( r_1 , \ldots , r_n ; \ee (1), \ldots ,
\ee (n) ) }
\newcommand{\ppree}{ {\cal P} ( r_1 , \ldots , r_{2p}; \ee (1), \ldots ,
\ee (2p) ) }
\newcommand{\qttheta}{ {\cal Q} ( t_{1}, \ldots , t_{n}; \theta_{1}, 
\ldots , \theta_{n} ) }
\newcommand{\qottheta}{ {\cal Q}( t_{1}, \ldots , t_{n}; \theta_{1}, 
\ldots , \theta_{n} ) }
\newcommand{\etatpi}{ \eta ( t_{1}, \ldots , t_{n} ; \pi ) }
\newcommand{\rhoijs}{ ( \rho_{i,j} )_{i,j} }
\newcommand{\wns}{ ( U_{n} )_{n=1}^{\infty} }

\newcommand{\downrightarrow}{\rightarrow}
\newcommand{\uprightarrow}{\rightarrow}

\newcommand{\A}{ {\cal A} }
\newcommand{\BB}{ {\cal B} ( \Omega , \A ) }
\newcommand{\C}{ \mbox{\bf C} }
\newcommand{\M}{ {\cal M} }
\newcommand{\N}{ {\cal N} }
\newcommand{\R}{ \mbox{\bf R} }
\newcommand{\T}{ {\cal T} }
\newcommand{\tq}{ {\cal T}_{q} }
\newcommand{\Z}{ \mbox{\bf Z} }

\newcommand{\ps}{ ( \Omega , {\cal F} , P ) }
\newcommand{\ncps}{ ( {\cal A} , \varphi ) }
\newcommand{\ncpsk}{ ( {\cal A}_{k} , \varphi_{k} ) }

\begin{document}

\title{\bf Random unitaries in non-commutative tori, \\ 
and an asymptotic model for q-circular systems}

\author{
James A. Mingo $^{*}$ \\
Department of Mathematics and Statistics  \\
Queen's University  \\
Kingston, Ontario K7L 3N6, Canada \\
mingoj@mast.queensu.ca
\and
Alexandru Nica 
\thanks{ Research supported by a grant from the 
Natural Sciences and Engineering Research Council, Canada.} \\
Department of Pure Mathematics \\
University of Waterloo \\
Waterloo, Ontario N2L 3G1, Canada \\
anica@math.uwaterloo.ca}

\date{ }

\maketitle

\begin{abstract}
We consider the concept of {\em q-circular system}, which is a 
deformation of the circular system from free probability, taking place
in the framework of the so-called ``$q$-commutation relations''.
We show that certain averages of random unitaries in non-commutative 
tori behave asymptotically like a $q$-circular system. More precisely:
let $q$ be in $(-1,1)$; let $s,k$ be positive integers; let 
$( \rho_{ij} )_{1 \leq i<j \leq ks}$ be independent random variables
with values in the unit circle, such that $\int \rho_{ij} = q$,
$\forall \ 1 \leq i<j \leq ks$; and let $U_{1}, \ldots , U_{ks}$ be 
random unitaries such that $U_{i}U_{j} = \rho_{ij} U_{j}U_{i}$, 
$\forall 1 \leq i<j \leq ks$. If we set:
\[
X_{r} \ := \ \frac{1}{\sqrt k} ( \ U_{r} + U_{r+s} + \cdots +
U_{r+(k-1)s} \ ), \ \ 1 \leq r \leq s,
\]
then the family $X_{1}, \ldots , X_{s}$ behaves for 
$k \rightarrow \infty$ like a $q$-circular system with $s$ elements.

The above result generalizes to the case when instead of the hypothesis
``$\int \rho_{ij} = q$'' we start with ``$\int \rho_{ij} = z$'', 
where $z$ is a complex number such that $|z| < 1$. In this case the 
limit distribution of $X_{1}, \ldots , X_{s}$ is what we call a 
{\em  z-circular system}. From the combinatorial point of view, the 
new feature brought in by a $z$-circular system is that its description
involves the enumeration of {\em oriented} crossings of certain pairings;
it is only in the case when $z= \overline{z} =q$ that the orientations 
cancel out, allowing the $q$-circular system to be described via 
non-oriented crossings.

As a consequence of the result, one can easily construct families of 
random matrices which converge in distribution to $q$-circular (or more 
generally $z$-circular) systems.
\end{abstract}

\newpage

\setcounter{section}{1}
\setcounter{equation}{0}
{\large\bf 1. Introduction and statement of the results} 

$\ $

This section is divided into subsections as follows:

\vspace{10pt}

1.1 The concept of circular system.

1.2 $q$-circular systems.

1.3 An asymptotic model for $q$-circular systems.

1.4 Oriented crossings.

1.5 $z$-circular systems.

1.6 Refinements of Theorem 1.5.3.

1.7 Approximation with random matrices.

$\ $

$\ $

{\bf 1.1 The concept of circular system} was introduced by D. Voiculescu in
\cite{V2}, and plays an important role in his theory of free probability.
The definition goes as follows. Let $\ncps$ be a 
{\em $C^{*}$-probability space} -- by which we mean that $\A$ is a unital
$C^{*}$-algebra and $\varphi$ is a state of $\A$ ($\varphi : \A \rightarrow 
\C$ positive linear functional, such that $\varphi (I) = 1$). The elements 
$c_{1}, \ldots , c_{s} \in \A$ ($s \geq 1$) are said to form a 
{\em circular system} in $\ncps$ if the family
\begin{equation}
\frac{c_{1}+c_{1}^{*}}{\sqrt{2}}, 
\frac{c_{1}-c_{1}^{*}}{i \sqrt{2}},  \ldots ,
\frac{c_{s}+c_{s}^{*}}{\sqrt{2}}, 
\frac{c_{s}-c_{s}^{*}}{i \sqrt{2}}
\end{equation}
is {\em free} in $\ncps$, and if each of the selfadjoint elements listed in
(1.1) has normalized semicircular distribution with respect to $\varphi$.
The fact that an element $a=a^{*} \in \A$ has normalized semicircular 
distribution means by definition that 
\begin{equation}
\varphi (a^{n} ) \ = \ \frac{1}{2 \pi} \int_{-2}^{2} t^{n} \sqrt{4-t^{2}}
\ dt, \ \forall n \geq 1.
\end{equation}
For the definition of freeness in $\ncps$, we refer to \cite{VDN}, Chapter 2.

The definition of a circular system given above can be rephrased in a 
purely combinatorial way, by indicating the general formula of the joint
moments of $c_{1},c_{1}^{*}, \ldots , c_{s},c_{s}^{*}$, i.e. of the 
expressions
\begin{equation}
\varphi ( \ c_{r_{1}}^{\ee (1)} \cdots c_{r_{n}}^{\ee (n)} \ ), \ 
n \geq 1, \ r_{1}, \ldots , r_{n} \in \{ 1, \ldots , s \} , \
\ee (1), \ldots , \ee (n) \in \{ 1, * \} .
\end{equation}
Namely, it turns out that every number in (1.3) is a non-negative integer,
which ``counts a certain family of non-crossing pairings''; this statement 
will be made precise (and generalized) in Section 1.2 below.

An important realization of a circular system, given in \cite{V2}, uses 
creation and annihilation operators on a full Fock space: if $\T$ is the 
full Fock space over $\C^{2s}$, if $\xi_{1}, \xi_{2}, \ldots , \xi_{2s}$ 
is an orthonormal basis of $\C^{2s}$, and if $l_{1},l_{2}, \ldots , l_{2s}$
are the creation operators on $\T$ determined by 
$\xi_{1}, \xi_{2}, \ldots , \xi_{2s}$, then
\begin{equation}
\frac{(l_{1}+l_{1}^{*})+i(l_{2}+l_{2}^{*})}{\sqrt{2}}, \ldots ,
\frac{(l_{2s-1}+l_{2s-1}^{*})+i(l_{2s}+l_{2s}^{*})}{\sqrt{2}} \ \in \
B( \T )
\end{equation}
form a circular system with respect to the so-called vacuum-state on
$B( \T )$. By performing an orthogonal transformation on
$\xi_{1}, \xi_{2}, \ldots , \xi_{2s}$ it is seen that instead of the family
in (1.4) one can also use 
\begin{equation}
l_{1}+l_{2}^{*}, \ldots , l_{2s-1}+l_{2s}^{*} \ \in \ B( \T ).
\end{equation}
The precise definitions of the objects involved in this realization of a 
circular system is reviewed in Section 2.1 below.

$\ $

$\ $

{\bf 1.2 q-circular systems.} In work related to $q$-deformations of the 
canonical commutation relations, Bo\.zejko and Speicher \cite{BS} obtained 
a remarkable deformation of the full Fock space, called the $q$-Fock space;
here $q$ is a parameter in $(-1,1)$, and the actual full Fock space is 
obtained for $q=0$ (see review in Section 2.2 below). In connection to 
this, Bo\.zejko and Speicher studied the
distribution -- called by them {\em q-Gaussian} -- of $l+l^{*}$, where $l$
is an appropriately normalized creation operator on the $q$-Fock space.
They discovered that the $q$-Gaussian distribution is the probability measure
associated to an important family of orthogonal polynomials, the 
$q$-continuous Hermite polynomials. If $q=0$, then the $q$-Gaussian is the 
semicircular distribution appearing in (1.2), while the usual Gaussian is 
obtained in the limit $q \rightarrow 1$.

From the work in \cite{BS} and its continuation in \cite{BKS} it is clear
that if in the formulas (1.4), (1.5) one makes $l_{1}, l_{1}^{*}, \ldots ,
l_{s}, l_{s}^{*}$ be creation/annihilation operators on the $q$-Fock 
space, then this should provide realizations of what one should call a 
$q$-circular system.

But how does one actually define a $q$-circular system? It is unfortunate
that the definition cannot be made in the same way as in the first 
paragraph of Section 1.1; this is because of the absence of a notion of
``$q$-freeness in $\ncps$''. However, the combinatorial reformulation 
mentioned in the second paragraph of Section 1.1 can be extended to the 
$q$-case. Indeed, it is possible to give an explicit combinatorial
description of the joint moments of the candidates of $q$-circular systems 
mentioned above (the families (1.4), (1.5), but where the $l_{i}$'s act
on the $q$-Fock space). We state this formally in the next definition and 
proposition.

Let us first succinctly describe the meaning of the combinatorial terms 
involved in Definition 1.2.2.

$\ $

{\bf 1.2.1 Notations.} The fact that $\pi = \{ B_{1}, \ldots , B_{p} \}$
is a {\em pairing} of $\{ 1, \ldots , n \}$ means that 
$B_{1} \cup \cdots \cup B_{p} = \{ 1, \ldots , n \}$, disjoint, and each
of $B_{1}, \ldots , B_{p}$ has exactly two elements. (Of course, $n$ must 
be even in order for $\{ 1, \ldots , n \}$ to have any pairings.) Two 
blocks $B_{i} = \{ a_{i}, b_{i} \}$ and $B_{j} = \{ a_{j}, b_{j} \}$ of 
a pairing $\pi = \{ B_{1}, \ldots , B_{p} \}$ are said to {\em cross} if
either $a_{i} < a_{j} < b_{i} < b_{j}$ or $a_{j} < a_{i} < b_{j} < b_{i}$;
the {\em number of crossings} of $\pi$ is 
\begin{equation}
cr ( \pi ) \ := \  \mbox{card} \{ (i,j) \ | \ 1 \leq i < j \leq p, \ 
B_{i} \mbox{ and } B_{j} \mbox{ cross} \} .
\end{equation}

$\ $

{\bf 1.2.2 Definition.} Let $\ncps$ be a $C^{*}$-probability space, and let 
$q$ be in $(-1, 1)$. The elements $c_{1}, \ldots , c_{s} \in \A$ ($s \geq 1$)
are said to form a {\em q-circular system} in $\ncps$ if for every 
$n \geq 1$, $\rs$, $\ees$ we have:
\begin{equation}
\varphi ( \cs ) \ = \ \sum_{\pi \in \pree} \ q^{cr ( \pi )} ,
\end{equation}
where $\pree$ denotes the set of all pairings $\pi$ =
$\{ \ \{ a_{1}, b_{1} \}, \ldots , \{ a_{p},b_{p} \} \ \}$ of 
$\{ 1, \ldots , n \}$ which have the property that $r_{a_{i}} = r_{b_{i}}$ 
and $\ee (a_{i}) \neq \ee (b_{i})$, $\forall \ 1 \leq i \leq p$. (In the 
case that ${\cal P} ( r_{1}, \ldots , r_{n};$ 
$\ee (1), \ldots , \ee (n) )$ is the empty set, the right-hand side 
of (1.7) is taken to be equal to 0.)

$\ $

{\bf 1.2.3 Proposition.} Let $q$ be in $(-1,1)$, and let $s$ be a positive
integer. Let $\tq$ denote the $q$-Fock space over $\C^{2s}$; consider an 
orthonormal basis $\xi_{1}, \ldots , \xi_{2s}$ of $\C^{2s}$ and let
$l_{1}, \ldots , l_{2s} \in B( \tq )$ be the creation operators associated
to $\xi_{1}, \ldots , \xi_{2s}$. Then the families of operators
\begin{equation}
\frac{(l_{1}+l_{1}^{*})+i(l_{2}+l_{2}^{*})}{\sqrt{2}}, \ldots ,
\frac{(l_{2s-1}+l_{2s-1}^{*})+i(l_{2s}+l_{2s}^{*})}{\sqrt{2}} \ \in \
B( \tq )
\end{equation}
and
\begin{equation}
l_{1}+l_{2}^{*}, \ldots , l_{2s-1}+l_{2s}^{*} \ \in \ B( \tq )
\end{equation}
are $q$-circular systems with respect to the vacuum-state on $B( \tq )$.

$\ $

Note that if $q=0$, then the right-hand side of (1.7) counts the non-crossing
pairings in $\pree$; this recaptures (and makes precise) the statement
following to Eqn.(1.3). In the particular case $q=0$, a proof of 
Proposition 1.2.3 can be made by using the concept of R-transform;
indeed, the R-transform of the family
$c_{1}, c_{1}^{*}, \ldots , c_{s}, c_{s}^{*}$ has a very simple form (see
e.g. Eqn.(1.6) of \cite{NS} for the case $s=1$), and the joint moments can be 
calculated from the knowledge of the R-transform.
For general $q \in (-1,1)$, the statement of Proposition 1.2.3 does not 
seem to have been previously considered, but can be inferred without 
difficulty from the results of \cite{BS} and \cite{BKS} (see Sections
2.3-2.5 below).

$\ $

$\ $

{\bf 1.3 An asymptotic model for q-circular systems.} We now arrive to the
main object of concern of the present paper,  which is a certain asymptotic
model for a $q$-circular system. The idea of using an asymptotic model for
a circular system (case $q=0$) was brought to fact by Voiculescu in
\cite{V1}, and then very successfully used in \cite{V2}. The asymptotic 
model observed in this paper is of a different nature than the one in
\cite{V1}, and is obtained by averaging unitaries in non-commutative tori.

$\ $

{\bf 1.3.1 Definition.} Let $q$ be in $(-1,1)$ and let $s$ be a positive 
integer. Suppose that for every $k \geq 1$ we are given a $C^{*}$-probability
space $\ncpsk$ and a family $c_{1;k}, \ldots, c_{s;k}$ of $\A_{k}$. 
We will say that these families {\em converge in distribution} to a 
$q$-circular system if for every $n \geq 1$, $\rs$, and $\ees$, the limit
\[
\lim_{k \rightarrow \infty} \ \varphi_{k} ( \
c_{r_{1};k}^{\ee (1)} \cdots c_{r_{n};k}^{\ee (n)} \ )
\]
exists and is equal to the right-hand side of Equation (1.7).

$\ $

{\bf 1.3.2 Definition.} Let $N$ be a positive integer, and let 
$( \rho_{i,j} )_{1 \leq i<j \leq N}$ be a family of complex numbers of 
absolute value 1. Let $\ncps$ be a $C^{*}$-probability space, and let 
$u_{1}, \ldots , u_{N}$ be elements of $\A$. We will say that 
$( u_{i} )_{i=1}^{N}$ is a $\rhoijs$--{\em commuting family of 
unitaries} if:

(i) every $u_{i}$ is a unitary, $1 \leq i \leq N$, and:

(ii) we have the relation $u_{i}u_{j} = \rho_{i,j} u_{j}u_{i}$, 
$\forall \ 1 \leq i<j \leq N$.

Moreover, we will say that 
$( u_{i} )_{i=1}^{N}$ is a $\rhoijs$--commuting {\em Haar} family of 
unitaries if in addition to (i) and (ii) we also have:

(iii) $\varphi ( u_{1}^{\lambda_{1}} \cdots
u_{N}^{\lambda_{N}} ) = 0$, $\forall \ ( \lambda_{1}, \ldots 
\lambda_{N} ) \in \Z^{N} \setminus \{ \ (0, \ldots , 0 ) \ \}$.

$\ $

$\rhoijs$--commuting Haar families of unitaries can be constructed for any 
choice of the $\rho_{i,j}$'s, and live naturally in a class of 
$C^{*}$-algebras called ``non-commutative tori'' -- see e.g. \cite{R}. 

$\ $

So now, let us fix a parameter $q \in (-1,1)$. We will denote:
\begin{equation}
\rho \ := \ q + i \sqrt{1-q^{2}} \ \ \ ( \ | \rho | = 1 \ ).
\end{equation}
We want to average families $u_{1}, \ldots , u_{N}$ of unitaries in a 
$C^{*}$-probability space, such that for every $1 \leq i<j \leq N$: either
$u_{i}$ and $u_{j}$ $\rho$-commute, or they $\rho^{-1}$-commute. Since
there is no canonical way to choose for which pairs $i<j$ we want to have
$u_{i}u_{j} = \rho u_{j}u_{i}$ and for which ones we want to have
$u_{i}u_{j} = \rho^{-1} u_{j}u_{i}$, we will use a ``randomization'' of
$u_{1}, \ldots , u_{N}$. That is, we will make $u_{1}, \ldots , u_{N}$
be random unitaries in a $C^{*}$-probability space, such that for every 
$1 \leq i<j \leq N$ we have:
\begin{equation}
P( \ u_{i}u_{j} = \rho u_{j}u_{i} \ ) \ = \ \frac{1}{2} \ = \ 
P( \ u_{i}u_{j} = \rho^{-1} u_{j}u_{i} \ ).
\end{equation}

This means that we will need the following version of Definition 1.3.2:

$\ $

{\bf 1.3.3 Definition.} Let $\ps$ be a probability space, let $N$ be a 
positive integer, and let $( \rho_{i,j} )_{1 \leq i<j \leq N}$ be a family
of random variables on $\Omega$ with values in $\{ \zeta \in \C \ |$
$| \zeta | =1 \}$. Let $\ncps$ be a $C^{*}$-probability space, where the 
$C^{*}$-algebra $\A$ is separable, and let $U_{1}, \ldots , U_{N}$ be 
measurable functions from $\Omega$ to $\A$. We will say that 
$( U_{i} )_{i=1}^{N}$ form a $\rhoijs$-{\em commuting family of 
random unitaries} in $\ncps$ if:

(j) $U_{i} ( \omega ) \in \A$ is a unitary, $\forall \ 1 \leq i \leq N$,
$\forall \ \omega \in \Omega$, and:

(jj) we have the relation $U_{i} ( \omega ) U_{j} ( \omega )$ = 
$\rho_{i,j} ( \omega ) U_{j} ( \omega ) U_{i} ( \omega )$, 
$\forall \ 1 \leq i<j \leq N$, $\forall \ \omega \in \Omega$.

Moreover, we will say that 
$( U_{i} )_{i=1}^{N}$ is a $\rhoijs$--commuting {\em Haar} family of 
random unitaries if in addition to (j) and (jj) we also have:

(jjj) $\varphi ( \ U_{1} ( \omega )^{\lambda_{1}} \cdots
U_{N} ( \omega )^{\lambda_{N}} \ ) = 0$, 
$\forall \ ( \lambda_{1}, \ldots \lambda_{N} )$ $\in$ 
$\Z^{N} \setminus \{ \ (0, \ldots , 0 ) \ \}$, 
$\forall \ \omega \in \Omega$.

$\ $

For the asymptotic model for a $q$-circular system it is sufficient to
consider $\rhoijs$--commuting families of random unitaries where the random 
variables $( \rho_{i,j} )_{1 \leq i<j \leq N}$ are independent, and each of 
them takes finitely many values. For such $\rho_{i,j}$'s, the interested 
reader should have no difficulty to verify that one can construct 
$\rhoijs$--commuting Haar families of random unitaries which live in a 
tensor product of non-commutative tori.

Since we will deal with random unitaries in a $C^{*}$-probability $\ncps$,
we will have to consider the new $C^{*}$-probability space where these 
random unitaries belong:

$\ $

{\bf 1.3.4 Notation.} Let $\ps$ be a probability space, and let $\ncps$ be a 
$C^{*}$-probability space, where the $C^{*}$-algebra $\A$ is separable.
We will denote by $\BB$ the set of all bounded measurable functions from 
$\Omega$ to $\A$. Then $\BB$ is a unital $C^{*}$-algebra, with the 
operations defined pointwise, and with the norm given by 
$|| f || := \sup \{ || f( \omega ) || \ | \ \omega \in \Omega \}$,
$f \in \BB$. Moreover, we have a natural state $E : \BB \rightarrow \C$
given by the formula
\begin{equation}
E(f) \ := \ \int_{\Omega} \varphi( \ f( \omega ) \ ) \ dP( \omega ), \ \ 
f \in \BB .
\end{equation}
It is immediate that $( \BB , E )$ is a $C^{*}$-probability space; also,
clearly, the unitaries in $\BB$ are random unitaries in $\A$, over the 
base space $\Omega$.
\footnote{ One could also consider the space $L^{\infty} ( \Omega , \A )$,
which is the quotient of $\BB$ by the relation of equality almost 
everywhere with respect to $P$. Since the estimates of moments done in
this paper are the same (no matter whether $\BB$ or 
$L^{\infty} ( \Omega , \A )$ is used), we prefer to stay with $\BB$. }

$\ $

We can now return to $q$ and $\rho$ of Equation (1.10), and state precisely
how an asymptotic model for the $q$-circular system is obtained.

$\ $

{\bf 1.3.5 Proposition.} Let $q$ be in $(-1,1)$, and let $s$ be a positive 
integer. Denote $\rho := q + i \sqrt{1-q^{2}}$. Suppose that for every 
$k \geq 1$ we have:

(a) A family $( \rho_{i,j;k} )_{1  \leq i<j \leq ks}$ of independent random
variables over some probability space $\Omega_{k}$, such that every 
$\rho_{i,j;k}$ takes only the values $\rho$ and $\rho^{-1}$, with 
$P( \ \rho_{i,j;k} = \rho \ ) = 1/2 = P( \ \rho_{i,j:k} = \rho^{-1} \ )$.

(b) A $( \rho_{i,j;k})_{i,j}$--commuting Haar family 
$U_{1;k}, \ldots , U_{ks;k}$ of random 
unitaries in some separable $C^{*}$-probability space $\ncpsk$.

Denote, for every $k \geq 1$:
\begin{equation}
X_{r;k} \ := \ \frac{1}{\sqrt{k}} ( U_{r;k} + U_{r+s;k} + \cdots +
U_{r+(k-1)s;k} ), \ \ 1 \leq r \leq s;
\end{equation}
then the family $( X_{1;k}, \ldots , X_{s;k} )$ converges in distribution, 
for $k \rightarrow \infty$, to a $q$-circular system.

$\ $

We should note here a similarity with the idea of the non-commutative 
central limit theorem of \cite{S}: in Proposition 1.3.5 we wrote $q$ as
a convex combination of $\rho$ and $\rho^{-1}$, whereas in \cite{S}
Speicher writes $q$ as a convex combination of $1$ and $-1$. In Theorem 
1.5.3 below we will generalize Proposition 1.3.5 to a case which contains 
both these situations, and where the only restriction on the 
$\rho_{i,j;k}$'s (besides their independence) concerns the values of 
their expectations. In order to state this more general result, we will
first introduce the concept of orientation for the crossings of a pairing.

$\ $

$\ $

{\bf 1.4 Oriented crossings.}

\vspace{10pt}

{\bf 1.4.1 Crossing of two segments.} We start from a simple geometric 
idea. Let $P,Q,U,V$ be distinct points in the plane, such that the 
segments $PQ$ and $UV$ cross. Consider the vector product 
$\vec{w} = \vec{PQ} \times \vec{UV}$, which is a vector perpendicular to
the plane of $P,Q,U,V$. If $\vec{w}$ is oriented upwards we will say that 
$PQ$ and $UV$ have a positive crossing, while if $\vec{w}$ is oriented
downwards we will say that $PQ$ and $UV$ have a negative crossing. In 
other words, if we denote the coordinates of $P$ by $(p_{1},p_{2})$,
the coordinates of $Q$ by $(q_{1},q_{2})$, etc, then the sign of the 
crossing between $PQ$ and $UV$ is equal to 
\begin{equation}
\mbox{sign} \Bigl( \ \mbox{det} \left(
\begin{array}{cc}
q_{1}-p_{1}  &  q_{2}-p_{2}  \\
v_{1}-u_{1}  &  v_{2}-u_{2}  
\end{array}  \right) \ \Bigr) .
\end{equation}

Note that the sign of the crossing is sensitive to the order of the points
of each segment, also to the order of the two segments; e.g, if $PQ$ and 
$UV$ have positive crossing then $QP$ and $UV$ have negative crossing, also
$UV$ and $PQ$ have negative crossing.

$\ $

{\bf 1.4.2 Crossings of a pairing.} Let now $n =2p$ be an even positive 
integer, and let $\pi = \{ B_{1}, \ldots , B_{p} \}$ be a pairing of 
$\{ 1, \ldots , n \}$. One can obtain a geometric representation of $\pi$,
by using the following recipe: draw a circle in the plane, and draw 
$n$ points $P_{1}, \ldots , P_{n}$ in counterclockwise order around the 
circle; then for every block $B_{i} = \{ a_{i}, b_{i} \}$ of $\pi$ draw the
line segment with endpoints $P_{a_{i}}$ and $P_{b_{i}}$. It is immediate 
that the blocks $B_{i}$ and $B_{j}$ cross if and only if the corresponding
line segments $P_{a_{i}}P_{b_{i}}$ and $P_{a_{j}}P_{b_{j}}$ do so. Thus 
the geometric representation of $\pi$ will display $p$ ( = $n/2$ )
line segments, which have a total number of $cr( \pi )$ points of 
intersection.

At this point, we would like to orient the crossings of the pairing $\pi$,
by using the considerations from 1.4.1. But in order to do so we need 
some additional data to be given, namely:

$( \alpha )$ a direction of running along the segment $P_{a_{i}}P_{b_{i}}$,
$1 \leq i \leq p$; and

$( \beta )$ an ordering of the $p$ segments $P_{a_{i}}P_{b_{i}}$,
$1 \leq i \leq p$.

It will be convenient to satisfy the above requirement $( \alpha )$
by giving a function $\ee : \{ 1, \ldots , n \} \rightarrow \{ 1,* \}$
with the property that for every block $B_{i} = \{ a_{i},b_{i} \}$ of
$\pi$ we have $\ee ( a_{i} ) \neq \ee ( b_{i} )$. In the presence of
such $\ee$, we will make the convention that every segment 
$P_{a_{i}}P_{b_{i}}$ is to be run from the point which is mapped by $\ee$
into $*$ towards the point which is mapped by $\ee$ into 1.

Concerning the requirement $( \beta )$, we will do the book-keeping by 
comparing the ordering of the blocks of $\pi$ which is used in the crossing
orientation against the ``standard'' ordering which lists the blocks in 
increasing order of their minimal elements. More precisely, let us assume
that the blocks $B_{1}, \ldots , B_{p}$ were from the beginning listed in 
standard order, with $\min (B_{1}) < \min (B_{2}) < \cdots < \min (B_{p})$.
Then giving an arbitrary ordering of the blocks amounts to giving a 
permutation $\sigma$ of the set $\{ 1, \ldots , p \}$: the convention 
we will use is that in the presence of such a permutation $\sigma$, the 
ordering ``$\prec$'' of the blocks of $\pi$ is defined such that
$B_{\sigma (1)} \prec B_{\sigma (2)} \prec \cdots \prec B_{\sigma (p)}$.

To summarize: we do not make the orientation of crossings for just the 
pairing $\pi$, but for a triple $( \orpi )$, where
$\ee : \{ 1, \ldots , n \} \rightarrow \{ 1,* \}$ has the property that 
$\ee ( a_{i} ) \neq \ee ( b_{i} )$ for every block 
$B_{i} = \{ a_{i},b_{i} \}$ of $\pi$, and $\sigma$ is a permutation of 
the set $\{ 1, \ldots , p \}$. For such $( \orpi )$, the orientation of 
crossings is achieved by drawing the geometric representation of $\pi$,
and then by using the method described in Section 1.4.1.

Let $( \orpi )$ be as in the preceding paragraph. We will denote by 
$cr_{+} ( \orpi )$ and $cr_{-} ( \orpi )$ the number of crossings of 
$( \orpi )$ which have positive, respectively negative, orientation. 
A distinctive feature of these numbers is of course that:
\begin{equation}
cr_{+} ( \orpi ) + cr_{-} ( \orpi ) \ = \ cr ( \pi ),
\end{equation}
the total number of crossings of $\pi$. We leave it as an exercise to the 
reader to check that if $\pi$ = $\{ B_{1}, \ldots , B_{p} \}$ with 
$B_{1}, \ldots , B_{p}$ listed in increasing order of their minimal 
elements, then the explicit formulas for $cr_{\pm} ( \orpi )$ are:
\[
cr_{+} ( \orpi ) \ = \ \mbox{card} \left\{ (i,j) \ 
\begin{array}{ll}
| & 1 \leq i<j \leq p, \ \mbox{$B_{i}$ and $B_{j}$ cross, }  \\
| & \ee ( \min (B_{i}) ) \cdot \ee ( \min (B_{j}) ) = \mbox{sign}
( \sigma (j) - \sigma (i) )
\end{array}
\right\}
\]
\begin{equation}
\end{equation}
\[
cr_{-} ( \orpi ) \ = \ \mbox{card} \left\{ (i,j) \ 
\begin{array}{ll}
| & 1 \leq i<j \leq p, \ \mbox{$B_{i}$ and $B_{j}$ cross, }  \\
| & \ee ( \min (B_{i}) ) \cdot \ee ( \min (B_{j}) ) = - \mbox{sign}
( \sigma (j) - \sigma (i) )
\end{array}
\right\} ,
\]
where in the products ``$\ee ( \min (B_{i})) \cdot \ee ( \min (B_{j}))$'' of
(1.16) the following convention is used: if we encounter a product 
of two symbols out of which at least one is a ``$*$'' (e.g. $* \cdot 1$,
or $* \cdot *$), then $*$ is to be treated like $-1$.

$\ $

$\ $

{\bf 1.5 z-circular systems.} The limit distribution which appears in the 
generalization of Proposition 1.3.5 is the following:

$\ $

{\bf 1.5.1 Definition.} Let $\ncps$ be a $C^{*}$-probability space, and 
let $z$ be a complex number such that $|z| < 1$. The elements 
$c_{1}, \ldots , c_{s} \in \A$ ($s \geq 1$) are said to form a 
{\em z-circular system} in $\ncps$ if:

-- for every positive odd integer $n$, for every $\rs$, and for every
$\ees$, we have that 
$\varphi ( c_{r_{1}}^{\ee (1)} \cdots c_{r_{n}}^{\ee (n)} ) = 0$; and

-- for every positive even integer $n =2p$, for every $\rs$, and for every
$\ees$, we have that:
\begin{equation}
\varphi ( \cs ) \ = \ \frac{1}{p!} \sum_{\sigma \in \symmp} \
\sum_{\pi \in \pree} \ z^{cr_{+} ( \orpi )} \cdot 
{\overline{z}}^{cr_{-} ( \orpi )} ,
\end{equation}
where $\symmp$ denotes the set of all permutations of $\{ 1, \ldots , p \}$,
the index set ${\cal P} ( r_{1}, \ldots , r_{n};$ 
$\ee (1), \ldots , \ee (n) )$ has the same meaning as in Definition 1.2.1,
and $cr_{\pm} ( \orpi )$ are as discussed in Section 1.4.2.

$\ $

{\bf 1.5.2 Remarks.} $1^{o}$ If $z=q \in (-1,1)$, then the concept of 
$z$-circular system reduces to the one of $q$-circular system from 
Definition 1.2.1. Indeed, in the relevant case of $n=2p$ appearing in 
Equation (1.17) we will now obtain:
\[
z^{cr_{+} ( \orpi )} \cdot {\overline{z}}^{cr_{-} ( \orpi )} \ = \ 
q^{cr_{+} ( \orpi ) + cr_{-} ( \orpi )} \ = \ q^{cr( \pi )},
\]
for every $\sigma \in \symmp$ and $\pi \in \pree$. So if one performs first
the summation over $\symmp$, then the right-hand side of Equation (1.17) 
reduces to the right-hand side of (1.7).

$2^{o}$ Starting from Definition 1.5.1, one can easily also define what
it means that a sequence of families 
$( c_{1;k}, \ldots , c_{s;k} )_{k \geq 1}$ converges in distribution to 
a $z$-circular system -- this is just an immediate adaptation of 
Definition 1.3.1.

$\ $

{\bf 1.5.3 Theorem.} Let $z$ be a complex number such that $|z|<1$,
and let $s$ be a positive integer. Suppose that for every 
$k \geq 1$ we have:

(a) A family $( \rho_{i,j;k} )_{1  \leq i<j \leq ks}$ of independent random
variables over some probability space $\Omega_{k}$, such that every 
$\rho_{i,j;k}$ takes values in the unit circle 
$\{ \zeta \in \C \ | \ | \zeta | =1 \}$, and has the property that
$ \int_{\Omega_{k}} \rho_{i,j;k} = z$.

(b) A $( \rho_{i,j;k} )_{i,j}$--commuting Haar family 
$U_{1;k}, \ldots , U_{ks;k}$ of random 
unitaries in some separable $C^{*}$-probability space $\ncpsk$.

Denote, for every $k \geq 1$:
\[
X_{r;k} \ := \ \frac{1}{\sqrt{k}} ( U_{r;k} + U_{r+s;k} + \cdots +
U_{r+(k-1)s;k} ), \ \ 1 \leq r \leq s;
\]
then the family $( X_{1;k}, \ldots , X_{s;k} )$ converges in distribution, 
for $k \rightarrow \infty$, to a $z$-circular system.

$\ $

{\bf 1.5.4 Remark.} 
$1^{o}$ In the case when the random variables $\rho_{i,j;k}$ of 
Theorem 1.5.3 take values in $\{ -1,1 \}$, we obtain a statement 
which is close to the framework of the central limit theorem of
\cite{S}. We note however that even in this case, the Proposition 1.6.2
below -- which generalizes Theorem 1.5.3, and is the statement that 
we really prove -- does not follow from the results of \cite{S}.

$2^{o}$ There are some natural questions which are raised by
the preceding theorem, concerning the possibility of realizing a 
$z$-circular system as a family of operators on some Hilbert space.

One approach that can be used is the following. Consider the unital algebra
$\C \langle X_{1},Y_{1}, \ldots ,$
$X_{s},Y_{s} \rangle =: \A_{o}$ of
polynomials in $2s$ non-commuting indeterminates 
$X_{1},Y_{1}, \ldots , X_{s},Y_{s}$, and make $\A_{o}$ be a $*$-algebra 
by introducing on it the (uniquely determined) $*$-operation with the 
property that $X_{r}^{*} = Y_{r}$, $1 \leq r \leq s$. Let us moreover 
consider the linear functional $\varphi_{o} : \A_{o} \rightarrow \C$
determined by the fact that $\varphi_{o} (1) = 1$ and that
$\varphi_{o} ( X_{r_{1}}^{\ee (1)} \cdots X_{r_{n}}^{\ee (n)} )$ is equal 
to the right-hand side of Equation (1.17), for every 
$n \geq 1$, $\rs$, $\ees$.
The Theorem 1.5.3 ensures that $\varphi_{o}$ is a positive functional
$( \varphi_{o} ( P^{*}P ) \geq 0, \ \forall \ P \in \A_{o} )$; indeed, 
it is easy to rephrase the theorem in a way which presents $\varphi_{o}$
as a pointwise limit of linear functionals $( \varphi_{k} )_{k=1}^{\infty}$,
each of the $\varphi_{k}$'s being positive.
But then one can consider the GNS construction for $\varphi_{o}$; this should
yield a $*$-representation $\Phi_{o} : \A_{o} \rightarrow B( {\cal H} )$,
with a cyclic vector $\xi_{o} \in {\cal H}$, such that
$\Phi_{o} (X_{1}), \ldots , \Phi_{o} (X_{s})$ form a $z$-circular system 
with respect to the vector-state on $B( {\cal H} )$ given by $\xi_{o}$.
The point we cannot settle here is whether the operators 
$\Phi_{o} (X_{1}), \ldots , \Phi_{o} (X_{s})$ are indeed bounded on 
${\cal H}$. We believe nevertheless that this is true, and that the 
condition ``$|z|<1$'' from the definition of a $z$-circular system should
be essential in proving it.

Another approach which can be tried in order to realize $z$-circular systems
would be by generalizing the Proposition 1.2.3 to the framework of an
appropriately defined $z$-Fock space. The concept of $q$-Fock space which 
will be reviewed in Section 2.2 below was amply generalized in 
\cite{BS2}, \cite{JSW}; on the other hand, Fock space constructions 
related to the framework of \cite{S} (spin systems with mixed commutation
and anti-commutation relations) are discussed in \cite{B}. It isn't however
clear if any of these constructions can be tailored to give a $z$-Fock space 
as required by the situation at hand.

$\ $

$\ $

{\bf 1.6 Refinements of Theorem 1.5.3.} It is useful (for instance for the 
approximation with random matrices shown in Section 1.7 below) to note that
one can relax some of the hypotheses of Theorem 1.5.3, and still obtain the 
same conclusion. In the next theorem we weaken the hypotheses on the 
expectations $\int \rho_{i,j}$, and on the Haar condition (jjj) from 
Definition 1.3.3. The weakened Haar condition is described as follows: let
$( U_{i} )_{i=1}^{N}$ be a $( \rho_{i,j} )_{i,j}$--commuting family of random
unitaries, in the sense of (j)$+$(jj) of Definition 1.3.3, and let $L$ be a
positive integer. We will say that $( U_{i} )_{i=1}^{N}$ is an 
{\em L-mimic of a Haar family} if it satisfies:
\[
\mbox{(jjj-$L$)} \ \ 
\varphi ( \  U_{1} ( \omega )^{\lambda_{1}} \cdots 
U_{N} ( \omega )^{\lambda_{N}}  \ ) \ = \ 0, \ \ \forall \ 
( \lambda_{1}, \ldots , \lambda_{N} ) \in ( \ (-L,L) \cap \Z \ )^{N}
\setminus \{ (0, \ldots , 0) \} .
\]

$\ $

{\bf 1.6.1 Proposition.} Let $z$ be a complex number such that $|z|<1$,
and let $s$ be a positive integer. Let $( \delta_{k} )_{k=1}^{\infty}$ be 
a sequence of positive real numbers, and let $( L_{k} )_{k=1}^{\infty}$ be 
a sequence of positive integers, such that $\delta_{k} \downrightarrow 0$
and $L_{k} \uprightarrow \infty$. Suppose that for every $k \geq 1$ we have:

(a) A family $( \rho_{i,j;k} )_{1  \leq i<j \leq ks}$ of independent 
random variables over some probability space $\Omega_{k}$, such that 
every $\rho_{i,j;k}$ takes values in the unit circle and has the property 
that $| z - \int_{\Omega_{k}} \rho_{i,j;k} | \leq \delta_{k}$.

(b) A $( \rho_{i,j;k} )_{i,j}$--commuting family $U_{1;k}, \ldots , U_{ks;k}$ 
of random unitaries in some separable $C^{*}$-probability space $\ncpsk$,
such that $U_{1;k}, \ldots , U_{ks;k}$ is an $L_{k}$-mimic of a Haar family.

Denote, for every $k \geq 1$:
\[
X_{r;k} \ := \ \frac{1}{\sqrt{k}} ( U_{r;k} + U_{r+s;k} + \cdots +
U_{r+(k-1)s;k} ), \ \ 1 \leq r \leq s;
\]
then the family $( X_{1;k}, \ldots , X_{s;k} )$ converges in distribution, 
for $k \rightarrow \infty$, to a $z$-circular system.

$\ $

It is worth recording that the statement of 1.6.1 follows from an estimate 
of moments which can be formulated simply, as described in the next 
proposition. (Since the extra indices ``$k$'' are not necessary in
Proposition 1.6.2, we will write 
in its statement $\rho_{i,j}, U_{i},X_{r}$ instead of 
$\rho_{i,j;k},U_{i;k},X_{r;k}$, respectively.)

$\ $

{\bf 1.6.2 Proposition.} Let $z$ be a complex number such that $|z|<1$,
and let $s$ be a positive integer. Let $\delta$ be a positive real number, 
and let $L$ be a positive integer. Let $k$ be a positive integer, and 
suppose that we have:

(a) A family $( \rho_{i,j} )_{1  \leq i<j \leq ks}$ of independent random
variables over some probability space $\Omega$, such that every $\rho_{i,j}$ 
takes values in the unit circle and has the property that
$| z - \int_{\Omega} \rho_{i,j} | \leq \delta$.

(b) A $( \rho_{i,j} )_{i,j}$--commuting family $U_{1}, \ldots , U_{ks}$ 
of random unitaries in some separable $C^{*}$-probability space $\ncps$,
such that $U_{1}, \ldots , U_{ks}$ is an $L$-mimic of a Haar family.

We denote:
\begin{equation}
X_{r} \ := \ \frac{1}{\sqrt{k}} ( U_{r} + U_{r+s} + \cdots +
U_{r+(k-1)s} ), \ \ 1 \leq r \leq s.
\end{equation}
We denote by $E : \BB \rightarrow \C$ the linear functional defined as in
Equation (1.12) of Notation 1.3.4. Then:

\vspace{10pt}

$1^{o}$ For every odd positive integer $n < L$, for every $\rs$, and for
every $\ees$, we have that:
\begin{equation}
E( \ X_{r_1}^{\ee (1)} \cdots X_{r_n}^{\ee (n)} \ ) \ = \ 0.
\end{equation}

\vspace{10pt}

$2^{o}$ For every even positive integer $n=2p$ such that
$n < \min (L,2k)$, for every $\rs$, and for every $\ees$, we have that:
\[
| \ E( \ X_{r_1}^{\ee (1)} \cdots X_{r_n}^{\ee (n)} \ ) \ - \ 
\frac{1}{p!} \sum_{\sigma \in \symmp} \
\sum_{\pi \in {\cal P} ( r_{1}, \ldots , r_{n}; \ee (1), \ldots , \ee (n) )}
\ z^{cr_{+} ( \orpi )} \cdot 
{\overline{z}}^{cr_{-} ( \orpi )}  \ | 
\]
\begin{equation}
< \ (2p+1)! \cdot ( \frac{1}{k} + \delta ).
\end{equation}

$\ $

The framework of Theorem 1.5.3 contains in particular the situation when the 
families of random unitaries $( U_{1;k}, \ldots , U_{ks;k} )_{k \geq 1}$
extend each other, i.e. when $U_{i;k+1} = U_{i;k}$, $\forall \ k \geq 1$,
$\forall \ 1 \leq i \leq ks$. When moving to the more general framework of
1.6.1, the case of the extending families of random unitaries needs to be 
discussed separately. One possibility of treating this case is provided 
by the following proposition.

$\ $

{\bf 1.6.3 Proposition.} Let $z$ be a complex number such that $|z|<1$, 
and let $s$ be a positive integer. Suppose that we have a family 
$( \rho_{m,n} )_{1 \leq m<n}$ of independent random variables with values 
in the unit circle, and a family $( U_{n} )_{n=1}^{\infty}$ of random 
unitaries in a separable $C^{*}$-probability space $\ncps$ (all the 
$\rho_{m,n}$'s and $U_{n}$'s defined on the same probability space 
$\Omega$), such that the following conditions are satisfied.

(a) The commutation relation
\begin{equation}
U_{m} ( \omega ) U_{n} ( \omega ) \ = \ 
\rho_{m,n} ( \omega ) U_{n} ( \omega ) U_{m} ( \omega )
\end{equation}
holds for every $1 \leq m <n$ and for every $\omega \in \Omega$.

(b) For every $\delta > 0$ there exists $m_{o} \geq 1$ such that:
$m_{o} \leq m <n$ $\Rightarrow$ 
$| z - \int_{\Omega} \rho_{m,n} | \leq \delta$.

(c) For every positive integer $L$ there exists $m_{o} \geq 1$ such that:
$m_{o} \leq m <n$ $\Rightarrow$ the family $U_{m}, \ldots , U_{n}$ is an 
$L$-mimic of a Haar family.

(d) If $n \geq 1$, $\lambda_{1}, \ldots , \lambda_{n} \in \Z$, 
$\omega \in \Omega$, and if at least one of 
$\lambda_{1}, \ldots , \lambda_{n}$ is equal to $\pm 1$ or to $\pm 2$,
then $\varphi ( \ U_{1} ( \omega )^{\lambda_{1}} \cdots
U_{n} ( \omega )^{\lambda_{n}} \ )$ = 0.

\vspace{10pt}

For every $k \geq 1$ we denote:
\begin{equation}
X_{r;k} \ = \ \frac{1}{\sqrt{k}} ( \ U_{r} + U_{r+s} + \cdots +
U_{r+(k-1)s} \ ), \ \ 1 \leq r \leq s;
\end{equation}
then the family $( X_{1;k}, \ldots , X_{s;k} )$ converges in distribution,
for $k \rightarrow \infty$, to a $z$-circular system.

$\ $

$\ $

{\bf 1.7 Approximation with random matrices.} We will now point out that, as 
a consequence of the results presented in Section 1.6, one can easily obtain 
families of 
{\em random matrices} which converge in distribution to a $z$-circular system.
In fact, it is nice to realize all these random matrices as random elements
in the same $C^{*}$-algebra, which will be an UHF-algebra (i.e. a certain
inductive limit of matrix algebras). 

So, let us fix a complex number $z$ such that $|z| < 1$. There exist unique
$\rho$ and $\gamma$ such that $| \rho | =1$, Im$(\rho) >0$, 
$\gamma \in (0,1)$, and $z = \gamma \rho + (1 - \gamma ) \overline{\rho}$.
Let us also fix a sequence $( \theta_{n} )_{n=1}^{\infty}$ of rational 
numbers in $(0,1)$, such that 
$\lim_{n \rightarrow \infty} e^{2 \pi i \theta_{n}} = \rho$, and such that
when we write $\theta_{n} = a_{n}/b_{n}$ with $a_{n},b_{n}$ relatively prime
positive integers, we get that $3 \leq b_{1} < b_{2} < \cdots$
$ < b_{n} < \cdots$.

For every $n \geq 1$, let us consider the finite dimensional $C^{*}$-algebra
\[
\A_{n} \ := \ M_{b_{1}} ( \C ) \otimes \cdots \otimes M_{b_{n}} ( \C ) ,
\]
(where $b_{j}$ is the denominator of $\theta_{j}$, as above); on 
$\A_{n}$ we consider the state $\varphi_{n}$ which is the tensor product of 
the normalized trace-functionals on 
$M_{b_{1}} ( \C ), \ldots , M_{b_{n}} ( \C )$. Let us furthermore consider 
the inductive limit 
\[
\A \ := \ \lim_{n \rightarrow \infty} \A_{n} , 
\]
where the mapping from $\A_{n-1}$ to $\A_{n}$ is 
$x \mapsto x \otimes 1_{b_{n}}$, for every $n \geq 1$. For notational 
convenience we shall regard each $\A_{n}$ as a unital subalgebra of 
$\A$. We will denote by $\varphi$ the unique state of $\A$ with the 
property that $\varphi | \A_{n} = \varphi_{n}$, $\forall \ n \geq 1$.

On the other hand, let us consider a probability space $\ps$ on which an
infinite family $( \xi_{m,n} )_{1 \leq m<n}$ of independent random variables
is given, such that every $\xi_{m,n}$ takes only the values $\pm 1$, with 
$P( \ \xi_{m,n} =1 \ ) \ = \ \gamma$, 
$P( \ \xi_{m,n} =-1 \ ) \ = \ 1- \gamma$. We denote
\begin{equation}
\rho_{m,n} ( \omega ) \ = \ e^{2 \pi i \theta_{m} \xi_{m,n} ( \omega )},
\ \ \forall \ \omega \in \Omega , \ \forall \ 1 \leq m <n.
\end{equation}
Then $( \rho_{m,n} )_{1 \leq m<n}$ is also an independent family of random
variables on $\Omega$, with values in the unit circle.

We construct a sequence of random unitaries $\wns$ in the $C^{*}$-algebra
$\A$, as follows. For every $n \geq 1$, consider first the 
$b_{n} \times b_{n}$-matrices:
\begin{equation}
V_{n} =  \left(  \begin{array}{ccccc}
0  &     &        &        &  1  \\
1  &  0  &        &        &     \\
   &  1  &  0     &        &     \\
   &     & \ddots & \ddots &     \\
   &     &        &   1    &  0  
\end{array} \right) , \ \ 
W_{n} = \left(  \begin{array}{ccccc}
e^{2 \pi i \theta_{n}} &                        &        &        &   \\
                       & e^{4 \pi i \theta_{n}} &        &        &   \\
                       &                        & \ddots &        &   \\
               &          &     & e^{2(b_{n}-1) \pi i \theta_{n}} &   \\
               &          &     &                                 & 1 
\end{array}  \right) ;
\end{equation}
then set for every $\omega \in \Omega$:
\begin{equation}
\left\{  \begin{array}{lll}
U_{1} ( \omega ) & = & V_{1}   \\
U_{2} ( \omega ) & = & W_{1}^{\xi_{1,2} ( \omega )} \otimes V_{2} \\
                 & \vdots &                                       \\
U_{n} ( \omega ) & = & W_{1}^{\xi_{1,n} ( \omega )} \otimes 
W_{2}^{\xi_{2,n} ( \omega )} \otimes \cdots \otimes
W_{n-1}^{\xi_{n-1,n} ( \omega )} \otimes V_{n}                    \\
                 & \vdots &                                       \\
\end{array}  \right.
\end{equation}
Clearly, the random unitary $U_{n}$ takes values in the finite dimensional
subalgebra $\A_{n} \subset \A$ (hence it is in fact a random unitary matrix
of size $b_{1}b_{2} \cdots b_{n}$).

We claim that:

$\ $

{\bf 1.7.1 Proposition.} The random variables $( \rho_{m,n} )_{1 \leq m <n}$
defined in Equation (1.23), and the random unitaries $(U_{n})_{n=1}^{\infty}$
defined in (1.25), satisfy the conditions considered in the Proposition 1.6.3.

$\ $

Indeed, the commutation relations (1.21) follow immediately from the fact 
that the matrices $V_{n},W_{n}$ in (1.24) satisfy the relation
$V_{n}W_{n} = e^{2 \pi i \theta_{n}} W_{n}V_{n}$. It is also 
immediate that:
\[
\Bigl| \ z -  \int_{\Omega} \rho_{m,n} \ \Bigr| \ = \  
| \ ( \ \gamma \rho + (1- \gamma ) \overline{\rho} \ ) -
( \gamma e^{2 \pi i \theta_{m}} +
(1- \gamma ) e^{-2 \pi i \theta_{m}} \ ) \ | 
\]
\[
\leq \ | \ e^{2 \pi i \theta_{m}} - \rho \ |, \ \ \forall \ 1 \leq m<n,
\]
and this gives the required behavior for the expectations of the
$\rho_{m,n}$'s. Concerning the Haar conditions, we leave it as an exercise 
to the reader to check that for every $1 \leq m < n$, the family
$W_{m},W_{m+1}, \ldots , W_{n}$ is a $b_{m}$-mimic of a Haar family; the 
verification of both this statement and of the hypothesis (d) in 1.6.3 
reduce to the fact that matrices of the form $V_{n}^{\alpha}W_{n}^{\beta}$
with $0 \neq \alpha \in (-b_{n},b_{n}) \cap \Z$ and $\beta \in \Z$ have 
only zeros on the diagonal (and therefore have zero trace).

Hence the recipe presented in Equation (1.22) of Proposition 1.6.3 will 
lead to an asymptotic $z$-circular system living in $\ncps$ (and 
which consists in fact of random matrices with sizes tending to infinity).

Alternatively, one can fabricate random matrices which converge in 
distribution to a $z$-circular system by cutting out disjoint segments 
of the sequence $(U_{n})_{n=1}^{\infty}$, and by invoking the Proposition
1.6.1. For example, in order to produce an asymptotic $z$-circular 
system with $s=2$ elements, one can set:
\[
X_{1;1} = U_{1}, \ \ X_{2;1} = U_{2},
\]
\[
X_{1;2} = \frac{ U_{3}+U_{5} }{\sqrt{2}}, \ \ 
X_{2;2} = \frac{ U_{4}+U_{6} }{\sqrt{2}}, 
\]
\vspace{10pt}
\[
X_{1;3} = \frac{ U_{7}+U_{9} +U_{11} }{\sqrt{3}}, \ \ 
X_{2;3} = \frac{ U_{8}+U_{10}+U_{12} }{\sqrt{3}}, \ldots
\]
\vspace{10pt}
(in general, the construction of $X_{1;k}$ and $X_{2;k}$ will use the
segment $U_{k^{2}-k+1}, \ldots , U_{k^{2}+k}$ of the sequence $\wns$ ).

$\ $

$\ $

The rest of the paper is divided into two sections. In Section 2
we review the $q$-Fock space, and prove Proposition 1.2.3. In
Section 3 we present the estimates of moments which lead to the theorems
presented in the Sections 1.3-1.6.

$\ $

$\ $

$\ $

\setcounter{section}{2}
\setcounter{equation}{0}
{\large\bf 2. Combinatorics of the joint moments of 
q-creation/annihilation operators} 

$\ $

{\bf 2.1 Review of the full Fock space.} In this paper we use the 
full Fock space over $\C^{2s}$ ($s$ a fixed positive integer), which is:
\begin{equation}
\T \ := \ \C \oplus \Bigl( \ \oplus_{n=1}^{\infty} 
( \C^{2s} )^{\otimes n} \ \Bigr) 
\end{equation}
(orthogonal direct sum of Hilbert spaces). The number 1 in the first summand
$\C$ on the right-hand side of Eqn.(2.1) is called the {\em vacuum-vector},
and is denoted by $\Omega$. The vector-state determined by $\Omega$ on 
$B( \T )$ is called the {\em vacuum-state}, and will be denoted
by $\varphi_{vac}$ $( \ \varphi_{vac} (X)$ = 
$\langle X \Omega \ | \ \Omega \rangle$, $\forall X \in B( \T ) \ )$.

For every $\xi \in \C^{2s}$ we denote by $l( \xi )$ the 
{\em creation operator} determined by $\xi$ on $\T$, which is described by:
\begin{equation}
\left\{
\begin{array}{l}
l( \xi ) \Omega = \xi        \\
                             \\
\begin{array}{r}
l( \xi ) ( \eta_{1} \otimes \cdots \otimes \eta_{n} )  =  
\xi \otimes  \eta_{1} \otimes \cdots \otimes \eta_{n}  \\
\forall \ n \geq 1, \ \forall \ \eta_{1}, \ldots , \eta_{n} \in \C^{2s} .
\end{array}
\end{array}    \right.
\end{equation}
The adjoint of $l( \xi )$ is called the {\em annihilation operator}
determined by $\xi$, and acts by:
\begin{equation}
\left\{
\begin{array}{l}
l( \xi )^{*} \Omega = 0      \\
                             \\
\begin{array}{r}
l( \xi )^{*} ( \eta_{1} \otimes \cdots \otimes \eta_{n} ) =
\langle \eta_{1} \ | \ \xi \rangle  \ \eta_{2} \otimes \cdots \otimes 
\eta_{n}  \\
\forall \ n \geq 1, \ \forall \ \eta_{1}, \ldots , \eta_{n} \in \C^{2s} .
\end{array}
\end{array}    \right.
\end{equation}
It is immediately verified that we have:
\begin{equation}
l( \xi )^{*} l( \eta ) \ = \ \langle \eta \ | \ \xi \rangle  \ I, \ \ 
\xi , \eta \in \C^{2s}.
\end{equation}
It is occasionally convenient to fix an orthonormal basis 
$\xi_{1}, \ldots , \xi_{2s}$ of $\C^{2s}$, and denote $l_{i} := l( \xi_{i} )$,
$1 \leq i \leq 2s$. Then (2.4) gives us that 
\begin{equation}
l_{i}^{*} l_{j} \ = \ \delta_{ij} I, \ \ 1 \leq i,j \leq 2s,
\end{equation}
i.e. that $l_{1}, \ldots , l_{2s}$ form a family of Cuntz isometries 
(isometries with mutually orthogonal ranges).  It is such a family which 
was used in (1.4) and (1.5) of Section 1.1, presenting realizations of
a circular system.

$\ $

{\bf 2.2 Review of the q-Fock space.} Besides the positive integer $s$,
we now also fix a parameter $q \in (-1,1)$. For every $n \geq 1$ 
we introduce an inner product $\langle \cdot \ , \cdot \rangle_{q}$ on 
$( \C^{2s} )^{\otimes n}$, determined by the formula:
\begin{equation}
\langle \xi_{1} \otimes \cdots \otimes \xi_{n} , 
\eta_{1} \otimes \cdots \otimes \eta_{n} \rangle_{q} \ := \
\sum_{\sigma \in \symmn} q^{inv( \sigma )} \ 
\langle \xi_{1} \ | \ \eta_{\sigma (1)} \rangle \cdots
\langle \xi_{n} \ | \ \eta_{\sigma (n)} \rangle  ,
\end{equation}
for $\xi_{1}, \ldots , \xi_{n}, \eta_{1}, \ldots , \eta_{n} \in \C^{2s}$,
where $\symmn$ denotes the set of all permutations of $\{ 1, \ldots , n \}$,
and where $inv( \sigma )$ stands for the number of inversions of the 
permutation $\sigma$ ( $inv ( \sigma ) := | \ \{ (i,j) \ | \ 
1 \leq i<j \leq n, \ \sigma (i) > \sigma (j) \} \ |$ ). The fact that (2.6)
gives indeed an inner product was shown in \cite{BS}. We view 
$\langle \cdot \ , \cdot \rangle_{q}$ as a ``deformation'' of the usual
inner product on $( \C^{2s} )^{\otimes n}$ (which would correspond to the 
case when $q=0$).

The $q$-Fock space over $\C^{2s}$ is defined as 
\begin{equation}
\tq \ := \ \C \oplus \Bigl( \ \oplus_{n=1}^{\infty} 
( \ ( \C^{2s} )^{\otimes n} , \langle \cdot \  , \cdot \rangle_{q} \ ) \ \Bigr) 
\end{equation}
(orthogonal direct sum of Hilbert spaces). The vacuum-vector of $\tq$
and the vacuum-state on $B( \tq )$ are defined in exactly the same 
way as for the full Fock space (cf. Section 2.1).
For every $\xi \in \C^{2s}$ there exists a unique operator in $B( \tq )$, 
denoted by $l_{q} ( \xi )$, such that:
\begin{equation}
\left\{
\begin{array}{l}
l_{q} ( \xi ) \Omega = \xi        \\
                             \\
\begin{array}{r}
l_{q} ( \xi ) ( \eta_{1} \otimes \cdots \otimes \eta_{n} )  =  
\xi \otimes  \eta_{1} \otimes \cdots \otimes \eta_{n}  \\
\forall \ n \geq 1, \ \forall \ \eta_{1}, \ldots , \eta_{n} \in \C^{2s} ;
\end{array}
\end{array}    \right.
\end{equation}
its adjoint acts by the formulas:
\begin{equation}
\left\{
\begin{array}{l}
l_{q} ( \xi )^{*} \Omega = 0      \\
                             \\
\begin{array}{r}
l_{q} ( \xi )^{*} ( \eta_{1} \otimes \cdots \otimes \eta_{n} ) =
\sum_{m=1}^{n} q^{m-1} \langle \eta_{m} \ | \ \xi \rangle  
\eta_{1} \otimes \cdots \otimes \eta_{m-1} \otimes \eta_{m+1}
\otimes \cdots \otimes \eta_{n}  \\
\forall \ n \geq 1, \ \forall \ \eta_{1}, \ldots , \eta_{n} \in \C^{2s} .
\end{array}
\end{array}    \right.
\end{equation}
$l_{q} ( \xi )$ and $l_{q}( \xi )^{*}$ are called the {\em q-creation} 
and respectively the {\em q-annihilation} operator determined by $\xi$. 
Note that the formulas describing $l_{q} ( \xi )$ are identical to those
for $l( \xi )$ in Section 2.1, but that the situation is not the same 
concerning the adjoints. (This is possible because $\tq$ has an inner
product which is a deformation of the one on $\T$.)

Instead of (2.4), we now get that the $q$-creation and $q$-annihilation 
operators satisfy:
\begin{equation}
l_{q} ( \xi )^{*} l_{q} ( \eta ) \ = \ 
q l_{q} ( \eta ) l_{q} ( \xi )^{*}  + \langle \eta \ | \ \xi \rangle \ 
I, \ \ \xi , \eta \in \C^{2s};
\end{equation}
these are called ``the $q$-commutation relations''.
It is occasionally convenient to fix an orthonormal basis 
$\xi_{1}, \ldots , \xi_{2s}$ of $\C^{2s}$, and denote\
$l_{i} := l_{q} ( \xi_{i} )$, $1 \leq i \leq 2s$. The Eqn.(2.10) then gives 
us that 
\begin{equation}
l_{i}^{*} l_{j} \ = \ q l_{j} l_{i}^{*} + 
\delta_{ij} I, \ \ 1 \leq i,j \leq 2s.
\end{equation}
It is such a family of operators in $B( \tq )$ which was used in 
(1.8) and (1.9) of Section 1.2, presenting realizations of a 
$q$-circular system.

$\ $

We now turn to the proof of Proposition 1.2.3. We will stick to the 
framework introduced in Section 2.2, including an orthonormal basis 
$\xi_{1}, \ldots , \xi_{2s}$ of $\C^{2s}$ which is fixed until the end
of the Section  2, and for which we denote $l_{i} := l_{q} ( \xi_{i} )$, 
$1 \leq i \leq 2s$. The argument will rely on a combinatorial formula 
established in \cite{BS} for the joint moments of 
$l_{1},l_{1}^{*}, \ldots , l_{2s}, l_{2s}^{*}$ with respect to the
vacuum-state $\varphi_{vac}$ on $B( \tq )$. This formula is stated as
follows:

$\ $

{\bf 2.3 Proposition} {\em (cf. \cite{BS} Part I, Proposition 2 on page 
529).} For every $n \geq 1$, $\ts$, $\thetas$, we have that
\begin{equation}
\varphi_{vac} ( l_{t_{1}}^{\theta_{1}} l_{t_{2}}^{\theta_{2}}  \cdots 
l_{t_{n}}^{\theta_{n}}  ) \ = \
\sum_{\pi \in \qttheta} \ \ q^{cr( \pi )} ,
\end{equation}
where $\qttheta$ denotes the set of all pairings $\pi$ =
$\{ \ \{ a_{1}, b_{1} \}, \ldots , \{ a_{p},b_{p} \} \ \}$ of 
$\{ 1, \ldots , n \}$ which have the property that $t_{a_{i}} = t_{b_{i}}$, 
$\theta_{a_{i}} = *$, and $\theta_{b_{i}} = 1$, $\forall \ 1 \leq i \leq p$. 

$\ $

We will first discuss the family of 
elements appearing in the formula (1.9) of Proposition 1.2.3. 

$\ $

{\bf 2.4 Proposition.} If $c_{1} := l_{1}+l_{2}^{*}, \ldots ,
c_{s} := l_{2s-1}+l_{2s}^{*}$, then $c_{1}, \ldots , c_{s}$ is a 
$q$-circular system with respect to the vacuum-state on $B( \tq )$.

$\ $

{\bf Proof.} We fix $n \geq 1$, $\rs$, $\ees$ for which we will verify 
that Eqn.(1.7) holds.

In this proof it will be convenient to use the following notation: given
$\ts$, $\thetas$, we will write
\begin{equation}
\tthetas \ \prec \ \rees
\end{equation}
to mean that for every $1 \leq m \leq n$ the operator 
$l_{t_{m}}^{\theta_{m}}$ is one of the two terms which form 
$c_{r_{m}}^{\ee (m)}$. (For example: if $c_{r_{m}}^{\ee (m)}$ = 
$c_{3}^{*} =(l_{5}+l_{6}^{*})^{*}$, then $l_{t_{m}}^{\theta_{m}}$ has to
be either $l_{5}^{*}$ or $l_{6}$; i.e., if $r_{m} = 3$ and $\ee_{m} = *$, 
then it is part of (2.13) that we have either $t_{m}=5$ and $\theta_{m} =*$,
or $t_{m}=6$ and $\theta_{m} = 1$.)

It is clear that:
\[
\cs \ = \ \sum_{\begin{array}{r} 
{\scriptstyle \tthetas \prec} \\
{\scriptstyle \rees} 
\end{array}  } \ \ls ,
\]
hence
\[
\varphi_{vac} ( \cs ) \ = \ \sum_{\begin{array}{r} 
{\scriptstyle \tthetas \prec} \\
{\scriptstyle \rees} 
\end{array}  } \ \varphi_{vac} ( \ls )
\]
\begin{equation}
= \ \sum_{\begin{array}{r} 
{\scriptstyle \tthetas \prec} \\
{\scriptstyle \rees} 
\end{array}  } \ \  \sum_{\pi \in \qottheta} \ \ 
q^{cr( \pi )} \ \ \mbox{( by (2.12) ).}
\end{equation}

We will next prove that:
\begin{equation}
\cup_{ \begin{array}{r} 
{\scriptstyle \tthetas \prec} \\
{\scriptstyle \rees} 
\end{array}  } \ \ \qottheta \ = \ \pree , 
\end{equation}
disjoint union.

In order to verify (2.15), let us first observe that:
\begin{equation}
\qottheta \ \subset \ \pree ,
\end{equation}
whenever $\tthetas \prec \rees$. This is immediately seen by comparing the 
definition of $\pree$ (see Definition 1.2.2) with the one of $\qottheta$,
and by taking into account how 
``$\prec$'' works. The inclusion (2.16) gives the ``$\subset$'' part of
(2.15).

We now pass to ``$\supset$'' of (2.15). We pick a partition $\pi \in \pree$,
and we will construct $\ts$, $\thetas$ such that
\begin{equation}
\left\{  \begin{array}{l}
\tthetas  \prec \rees, \ \mbox{and}  \\
\pi \in \qottheta .
\end{array}  \right.
\end{equation}

Let $B = \{ a,b \}$, with $a<b$, be an arbitrary block of $\pi$. From the 
fact that $\pi \in \pree$, we get that $r_{a} = r_{b} =: r$, and 
$\ee (a) \neq \ee (b)$. If $\ee (a) = 1$ and $\ee (b) = *$, this means
that $c_{r_{a}}^{ \ee (a)} = c_{r} = l_{2r-1} + l_{2r}^{*}$,
$c_{r_{b}}^{ \ee (b)} = c_{r}^{*} = l_{2r-1}^{*} + l_{2r}$, and 
we choose: $t_{a} = t_{b} = 2r$, $\theta_{a} = *$, $\theta_{b} = 1$
(such that $l_{t_{a}}^{\theta_{a}}$ is a term of $c_{r}$, and 
$l_{t_{b}}^{\theta_{b}}$ is a term of $c_{r}^{*}$).
If $\ee (a) = *$ and $\ee (b) = 1$, this means
that $c_{r_{a}}^{ \ee (a)} = c_{r}^{*} = l_{2r-1}^{*} + l_{2r}$,
$c_{r_{b}}^{ \ee (b)} = c_{r} = l_{2r-1} + l_{2r}^{*}$, and 
we choose: $t_{a} = t_{b} = 2r-1$, $\theta_{a} = *$, $\theta_{b} = 1$.

When we make the choices for $t_{a},t_{b}, \theta_{a}, \theta_{b}$ as 
described in the preceding paragraph, and for every block of $\pi$, we obtain
some  $\ts$ and $\thetas$ such that (2.17) holds. This completes the 
proof of ``$\supset$'' in (2.15). It is also immediate (by inspecting again,
one by one, the blocks of the partition $\pi$ considered above) that the 
choices for $t_{1}, \ldots , t_{n}, \theta_{1} , \ldots , \theta_{n}$
such that $\tthetas \prec \rees$ and at the same time 
$\qottheta \ni \pi$ are uniquely determined; this proves the disjointness 
of the union in (2.15).

Finally, from (2.15) it follows that the expression in (2.14) is
\[
\sum_{\pi \in \pree} \ q^{cr( \pi )}, 
\]
which is exactly the desired expression for $\varphi_{vac} ( \cs )$.
{\bf QED}

$\ $

It only remains that we prove the $q$-circularity of the family appearing
in (1.8) of Proposition 1.2.3. By using arguments from \cite{BKS}, this 
can in fact be reduced to the $q$-circularity of (1.9), which was shown 
above.

$\ $

{\bf 2.5 Proposition.} If we denote:
\[
c_{1}' := \frac{(l_{1}+l_{1}^{*})+i(l_{2}+l_{2}^{*})}{\sqrt{2}}, \ldots ,
c_{s}' := \frac{(l_{2s-1}+l_{2s-1}^{*})+i(l_{2s}+l_{2s}^{*})}{\sqrt{2}},
\]
then $c_{1} ', \ldots , c_{s}'$ is a $q$-circular system with respect to 
the vacuum-state on $B( \tq )$.

$\ $

{\bf Proof.} Recall that ``$l_{k}$'' stands here for ``$l_{q} ( \xi_{k} )$'',
$1 \leq k \leq 2s$, where $\xi_{1}, \ldots , \xi_{2s}$ is an orthonormal 
basis of $\C^{2s}$ which was fixed prior to the Proposition 2.3. Consider
the vectors $\eta_{1} , \ldots , \eta_{2s} \in \C^{2s}$ defined by:
\begin{equation}
\eta_{2r-1}  \ = \ \frac{\xi_{2r-1} + \xi_{2r}}{\sqrt{2}} , \ \ 
\eta_{2r}    \ = \ \frac{\xi_{2r-1} - \xi_{2r}}{i \sqrt{2}} , \ \
1 \leq r \leq s,
\end{equation}
and let $T$ denote the unique operator in $B( \C^{2s} )$ such that 
$T \xi_{k} =  \eta_{k}$, $1 \leq k \leq 2s$. It is immediate that 
$\eta_{1} , \ldots , \eta_{2s}$ is an orthonormal basis of $\C^{2s}$,
hence that $T$ is an orthogonal transformation.

Let ${\cal V} \subset \C^{2s}$ be the real vector space spanned by 
$\xi_{1}, \ldots , \xi_{2s}$ (i.e. the set of vectors of the form 
$\sum_{k=1}^{2s} \lambda_{k} \xi_{k}$, with $\lambda_{1}, \ldots , 
\lambda_{2s} \in \R$), and let ${\cal W} \subset \C^{2s}$ be
the real vector space spanned by $\eta_{1}, \ldots , \eta_{2s}$. 
Moreover, let $\M , \N \subset B( \tq )$ denote the 
von Neumann algebras generated by 
$\{ l_{q} ( \xi ) + l_{q} ( \xi )^{*} \ | \ \xi \in {\cal V} \}$, 
and respectively by
$\{ l_{q} ( \eta ) + l_{q} ( \eta )^{*} \ | \ \eta \in {\cal W} \}$.
The Theorem 2.11 of \cite{BKS} gives us the existence 
of a unital $*$-homomorphism $\Phi : \M \rightarrow \N$, which
preserves the vacuum-state (i.e. 
$\varphi_{vac} ( \Phi (x)) = \varphi_{vac} (x)$, $\forall \ x \in \M$),
and such that:
\begin{equation}
\Phi ( \ l_{q} ( \xi ) + l_{q} ( \xi )^{*} \ ) \ = \ 
l_{q} ( T \xi ) + l_{q} ( T \xi )^{*}, \ \ \forall \ \xi \in {\cal V}.
\end{equation}

It is obvious that the operators $c_{1}', \ldots , c_{s}'$ defined in the
statement of the proposition belong to $\M$, and an immediate calculation
which uses (2.18), (2.19), and the linearity of $l_{q} ( \cdot )$ gives
that:
\begin{equation}
\Phi ( c_{r} ' ) \ = \ l_{2r-1} + l_{2r}^{*}, \ \ 1 \leq r \leq s.
\end{equation}
Denoting $c_{r} := l_{2r-1} + l_{2r}^{*}$, $1 \leq r \leq s$, we thus 
obtain that $c_{1}, \ldots , c_{s} \in \N$ and also (since $\Phi$ is a
$*$-homomorphism which preserves $\varphi_{vac}$) that:
\[
\varphi_{vac} ( \cs ) \ = \ 
\varphi_{vac} ( \ \Phi ( (c_{r_{1}} ')^{\ee (1)} \cdots 
( c_{r_{n}} ')^{\ee (n)} ) \ )
\]
\begin{equation}
= \ \varphi_{vac} ( \ (c_{r_{1}}')^{\ee (1)} \cdots 
( c_{r_{n}} ')^{\ee (n)} \ ),
\end{equation}
for every $n \geq 1$ and $\rs$, $\ees$. But then the conclusion of the 
current proposition follows from (2.21) and Proposition 2.4.
{\bf QED}

$\ $
 
$\ $
 
$\ $
 
\setcounter{section}{3}
\setcounter{equation}{0}
{\large\bf 3. Moment estimates leading to asymptotic z-circular systems}

$\ $

In this section we prove the results stated in the Sections 1.3-1.6 of the 
Introduction. It is clear that in fact only the Propositions 1.6.2 and 1.6.3
need to be proved (then 1.6.1, 1.5.3, 1.3.5 will follow). 

The bulk of the section will be devoted to the estimates of moments 
presented in Proposition 1.6.2. We fix, from this moment on and until the 
end of Section 3.6, the framework described in 1.6.2. We will first dispose 
of the easy case appearing in the part $1^{o}$ of the proposition.

$\ $

{\bf 3.1 Proof} {\em of part $1^{o}$ in Proposition 1.6.2.}
By substituting $X_{1}, \ldots , X_{s}$ from their
definition in Eqn.(1.18), then by expanding the sums 
and by using the definition of $E$, we get:
\[
E ( \ X_{r_{1}}^{\ee (1)} \cdots X_{r_{n}}^{\ee (n)}  \ ) \ = \
\]
\[
= \ \frac{1}{k^{n/2}} \cdot  \sum_{
\begin{array}{c}
{\scriptstyle 1 \leq i_{1}, \ldots , i_{n} \leq ks \ such \ that}  \\
{\scriptstyle i_{1}=r_{1}(mod \ s), \ldots , i_{n}=r_{n}(mod \ s) }
\end{array} }  \ \ 
E ( \ U_{i_{1}}^{\ee (1)} \cdots U_{i_{n}}^{\ee (n)} \ )
\]
\begin{equation}
= \ \frac{1}{k^{n/2}} \cdot \sum_{
\begin{array}{c}
{\scriptstyle 1 \leq i_{1}, \ldots , i_{n} \leq ks \ such \ that} \\
{\scriptstyle i_{1}=r_{1}(mod \ s), \ldots , i_{n}=r_{n}(mod \ s) }
\end{array} }  \ \ \int_{\Omega} \varphi ( \
U_{i_{1}} ( \omega )^{\ee (1)} \cdots 
U_{i_{n}} ( \omega )^{\ee (n)} \ ) \ dP( \omega ).
\end{equation}

We will show that:
\begin{equation}
\varphi ( \ 
U_{i_{1}} ( \omega )^{\ee (1)} \cdots 
U_{i_{n}} ( \omega )^{\ee (n)} \ ) \ = \ 0,
\end{equation}
for every $\omega \in \Omega$ and every 
$1 \leq i_{1}, \ldots , i_{n} \leq ks$ such that 
$i_{1} = r_{1} (mod \ s), \ldots , i_{n} = r_{n} (mod \ s)$. This and
(3.1) clearly imply the conclusion of the lemma.

So let us fix $\omega \in \Omega$ and 
$1 \leq i_{1}, \ldots , i_{n} \leq ks$ such that 
$i_{1} = r_{1} (mod \ s), \ldots , i_{n} = r_{n} (mod \ s)$. The 
commutation relations satisfied by $U_{1}, \ldots , U_{ks}$ (see 
condition (jj) in Definition 1.3.3) give us that
\begin{equation}
U_{i_{1}} ( \omega )^{\ee (1)} \cdots 
U_{i_{n}} ( \omega )^{\ee (n)} \ = \ c \ 
U_{1} ( \omega )^{\lambda_{1}} \cdots 
U_{ks} ( \omega )^{\lambda_{ks}} ,
\end{equation}
where $c$ is a constant of absolute value 1, and where 
$\lambda_{1}, \ldots , \lambda_{ks} \in [-n,n] \cap \Z \subset
(-L,L) \cap \Z$. It cannot be true that
$\lambda_{1} = \cdots = \lambda_{ks} = 0$, because:
\[
\lambda_{1} + \cdots + \lambda_{ks} \ = \ 
| \  \{ 1 \leq m \leq n \ | \ \ee (m) =1 \} \ | 
\ - \ | \ \{ 1 \leq m \leq n \ | \ \ee (m) =* \} \ | ,
\]
which is an odd number (indeed, $| \ \{ m \ | \ \ee (m) =1 \} \ |$ and 
$| \ \{ m \ | \ \ee (m) =* \} \ |$ must have different parities, since 
their sum is the odd number $n$). But then the condition (jjj-$L$)
introduced in Section 1.6 gives us that
$\varphi ( \ U_{1} ( \omega )^{\lambda_{1}} \cdots 
U_{ks} ( \omega )^{\lambda_{ks}}  \ )$ = 0, and (3.2) is obtained
by applying $\varphi$ to both sides of (3.3). {\bf QED}

$\ $

We now move towards the sensibly harder case 
discussed in part $2^{o}$ of Proposition 1.6.2. We will start by making a
number of preliminary considerations.

Unlike in the preceding proof, where we did not need to know what was 
the constant $c$ in Equation (3.3), the arguments in the sequel will 
require some information about such constants which arise from 
commutations. The next lemma will be used for that.

$\ $

{\bf 3.2 Lemma.} Let $p$ be a positive integer and let 
$\pi = \{ B_{1} , \ldots , B_{p} \}$ be a pairing of $\ppset$, where the 
blocks $B_{1}, \ldots , B_{p}$ of $\pi$ are listed in increasing order of 
their minimal elements. Let ${\cal C}$ be a unital algebra and let 
$V_{1}, \ldots , V_{p}$ be invertible elements of ${\cal C}$ which satisfy 
the commutation relations
\begin{equation}
V_{l}V_{m} \ = \ \gamma_{lm} V_{m}V_{l}, \ \ 1 \leq l < m \leq p,
\end{equation}
where the $\gamma_{lm}$'s are some complex numbers. Define 
$W_{1}, \ldots , W_{2p}$ according to the formula:
\begin{equation}
W_{i} \ = \ \left\{  \begin{array}{lccl}
V_{l}  &  \mbox{if} & i = \mbox{min} (B_{l}) & \mbox{(for some 
$1 \leq l \leq p)$}                                               \\
       &            &                        &                    \\
V_{m}^{-1}  &  \mbox{if} & i = \mbox{max} (B_{m}) & \mbox{(for some 
$1 \leq m \leq p$).}                                               \\
\end{array}  \right.
\end{equation}
Then we have
\begin{equation}
W_{1}W_{2} \cdots W_{2p} \ = \ \Bigl( \ \prod_{
\begin{array}{c}
{\scriptstyle 1 \leq l < m \leq p \ such} \\
{\scriptstyle that \ B_{l} \ crosses \ B_{m}} 
\end{array}  } \  \gamma_{lm} \  \Bigr)  \ I.
\end{equation}

$\ $

{\bf Proof.} By induction on $p$. The case $p=1$ is obvious
(both sides of (3.6) are equal to $I$). 

Let us assume the lemma true for $p-1$ and prove it for $p \geq 2$.
Let $\pi = \{ B_{1}, \ldots , B_{p} \}$, $V_{1}, \ldots , V_{p}$ and 
$W_{1}, \ldots , W_{2p}$ be as in the statement of the lemma. We write
explicitly $B_{p} = \{ a,b \}$, $a<b$ (recall that $B_{p}$ is the block of 
$\pi$ with the largest minimal element). Note that $\{ a+1, \ldots , b-1 \}$
coincides with the set of maximal elements of the blocks $B_{l}$
$(1 \leq l \leq p-1)$ which cross $B_{p}$. By using this observation,
the rule (3.5) for defining $W_{a+1}, \ldots , W_{b-1}$, and the 
commutation relations (3.4), we obtain that:
\begin{equation}
(W_{a+1} \cdots W_{b-1}) V_{p}^{-1} \ = \ 
\Bigl( \ \prod_{ \begin{array}{c}
{\scriptstyle 1 \leq l \leq p-1 \ such} \\
{\scriptstyle that \ B_{l} \ crosses \ B_{p}} 
\end{array}  } \ \gamma_{lp} \  \Bigr)  \
V_{p}^{-1} (W_{a+1} \cdots W_{b-1}).
\end{equation}

On the other hand, let us denote by $\pi_{o}$ the pairing which is obtained
from $\pi$ by deleting the block $B_{p}$ and by redenoting the elements of
$\{ 1, \ldots , 2p \} \setminus B_{p}$ as $1,2, \ldots , 2p-2$, in 
increasing order. The induction hypothesis applied to $\pi_{o}$ and 
$V_{1}, \ldots , V_{p-1}$ gives us that:
\begin{equation}
W_{1} \cdots W_{a-1} W_{a+1} \cdots W_{b-1} W_{b+1} \cdots W_{2p} \ = \ 
\Bigl( \ \prod_{
\begin{array}{c}
{\scriptstyle 1 \leq l < m \leq p-1 \ such} \\
{\scriptstyle that \ B_{l} \ crosses \ B_{m}} 
\end{array}  } \ \gamma_{lm} \  \Bigr)  \ I.
\end{equation}
But then:
\[
W_{1}W_{2} \cdots W_{2p} \ = \ 
(W_{1} \cdots W_{a-1}) V_{p} (W_{a+1} \cdots W_{b-1}) V_{p}^{-1} 
(W_{b+1} \cdots W_{2p}) 
\]
\[
= \ \Bigl( \ \prod_{ \begin{array}{c}
{\scriptstyle 1 \leq l \leq p-1 \ such} \\
{\scriptstyle that \ B_{l} \ crosses \ B_{p}} 
\end{array}  } \ \gamma_{lp} \  \Bigr)  \
(W_{1} \cdots W_{a-1}) V_{p} V_{p}^{-1} (W_{a+1} \cdots W_{b-1}) 
(W_{b+1} \cdots W_{2p}) 
\]
( by Equation (3.7) )
\[
= \ \Bigl( \ \prod_{
\begin{array}{c}
{\scriptstyle 1 \leq l < m \leq p \ such} \\
{\scriptstyle that \ B_{l} \ crosses \ B_{m}} 
\end{array}  } \ \gamma_{lm} \  \Bigr)  \ I
\ \ \mbox{( by Equation (3.8) ).}
\]
{\bf QED}

$\ $

In the estimates of moments which will be presented below, we will also 
use the following notation and lemma. The positive integers $p,s,k$
appearing in 3.3 and 3.4 are the ones given in the statement of Proposition
1.6.2.

$\ $

{\bf 3.3 Notation.} Let $j_{1}, \ldots , j_{p}$ be distinct numbers in 
$\{ 1, \ldots , ks \}$. We will denote by $ord( j_{1}, \ldots , j_{p} )$ 
the permutation $\sigma$ of $\pset$ which keeps track of the order of 
$j_{1}, \ldots , j_{p}$; that is, $\sigma$ is the 
unique bijection from $\pset$ to itself which has the property that 
\begin{equation}
\sigma (l) < \sigma (m) \ \ \Leftrightarrow  \ \ j_{l} < j_{m}, 
\ \ \forall \ l \neq m
\mbox{ in } \pset .
\end{equation}

$\ $

{\bf 3.4 Lemma.} Let $\sigma$ be a permutation of $\pset$, and let 
$t_{1}, \ldots , t_{p}$ be in $\{ 1, \ldots , s \}$. Consider the number:
\begin{equation}
N( \sigma ; t_{1}, \ldots , t_{p} ) \ = \ | \ \Bigl\{ \ 
( j_{1}, \ldots , j_{p} ) \ 
\begin{array}{ll}
|  &  1 \leq j_{1}, \ldots , j_{p} \leq ks, \\
|  &  ord( j_{1}, \ldots , j_{p} ) = \sigma , \\
|  &  j_{1}=t_{1} (mod \ s), \ldots , j_{p} = t_{p} (mod \ s)
\end{array}  \ \Bigr\} \ | .
\end{equation}
Then:
\begin{equation}
\left(  \begin{array}{cc} k  \\ p  \end{array}  \right) \ \leq
N( \sigma; t_{1}, \ldots , t_{p} ) \ \leq
\left(  \begin{array}{cc} k+p  \\ p  \end{array}  \right) .
\end{equation}

$\ $

{\bf Proof.} It is immediate that
\[
N( \sigma ; t_{1}, \ldots , t_{p} )\ = \ 
N( id ; t_{\sigma^{-1} (1)}, \ldots , t_{\sigma^{-1} (p)} ),
\]
where $id$ denotes the identity permutation. Due to this fact, it suffices
to verify (3.11) in the case when $\sigma = id$; i.e, it suffices to verify 
that for any choice of $t_{1}, \ldots , t_{p} \in \{ 1, \ldots , s \}$, the
set
\begin{equation}
{\cal S} \ := \ \Bigl\{ \ (j_{1}, \ldots , j_{p} ) \
\begin{array}{ll}
|  &  1 \leq j_{1}< \cdots < j_{p} \leq ks,   \\
|  &  j_{1}=t_{1} (mod \ s), \ldots , j_{p} = t_{p} (mod \ s)
\end{array}  \ \Bigr\} 
\end{equation}
has cardinality between 
$\left(  \begin{array}{cc} k  \\ p  \end{array}  \right)$ and
$\left(  \begin{array}{cc} k+p  \\ p  \end{array}  \right)$.

Let us denote $I_{1} = \{ 1, \ldots , s \}$,
$I_{2} = \{ s+1, \ldots , 2s \} , \ldots ,$ 
$I_{k} = \{ (k-1)s+1, \ldots , ks \}$. To every 
$(j_{1}, \ldots , j_{p})$ in the set ${\cal S}$ of (3.12) we can 
associate the $p$-tuple $(m_{1}, \ldots , m_{p})$, where 
$1 \leq m_{1} \leq m_{2} \leq \cdots \leq m_{p} \leq k$ are determined 
by the conditions:
\[
j_{1} \in I_{m_{1}}, j_{2} \in I_{m_{2}}, \ldots , 
j_{p} \in I_{m_{p}}.
\]
Then the map $(j_{1}, \ldots , j_{m}) \mapsto (m_{1}, \ldots , m_{p})$ 
is one-to-one; this is immediately implied by the fact that every 
$(j_{1}, \ldots , j_{p})$ in the set ${\cal S}$ of (3.12) has to satisfy
the conditions $j_{1} = t_{1} (mod \ s), \ldots , j_{p} = t_{p}(mod \ s)$.
We hence obtain that the cardinality of ${\cal S}$ is bounded above by
\[
| \ \Bigl\{ (m_{1}, \ldots , m_{p}) \ | \ 
1 \leq m_{1} \leq m_{2} \leq \cdots \leq m_{p} \leq k \Bigr\} \ | \ = \ 
\left(  \begin{array}{cc} k+p-1  \\ p  \end{array}  \right) \ \leq \
\left(  \begin{array}{cc} k+p  \\ p  \end{array}  \right) .
\]

On the other hand, the range of the map 
$(j_{1}, \ldots , j_{m}) \mapsto (m_{1}, \ldots , m_{p})$ considered in 
the preceding paragraph contains all the $p$-tuples 
$(m_{1}, \ldots , m_{p})$ with the property that 
$m_{1} < m_{2} < \cdots < m_{p}$. Indeed, if
$m_{1} < m_{2} < \cdots < m_{p}$, then there are unique 
$j_{1}, \ldots , j_{p} \in \{ 1, \ldots , ks \}$ such that: 
$j_{1} \in I_{m_{1}}$ and $j_{1} = t_{1} (mod \ s)$;
$j_{2} \in I_{m_{2}}$ and $j_{2} = t_{2} (mod \ s) ; \ldots ,$
$j_{p} \in I_{m_{p}}$ and $j_{p} = t_{p} (mod \ s)$. These 
$j_{1}, \ldots , j_{p}$ form an element of the set ${\cal S}$ of (3.12),
which is mapped to $(m_{1}, \ldots , m_{p})$. So we obtain that the 
cardinality of ${\cal S}$ is bounded below by:
\[
| \ \Bigl\{ (m_{1}, \ldots , m_{p}) \ | \ 
1 \leq m_{1} < m_{2} < \cdots < m_{p} \leq k \Bigr\} \ | \ = \ 
\left(  \begin{array}{cc} k  \\ p  \end{array}  \right) .
\]
{\bf QED}

$\ $

We are now ready to attack the proof of part $2^{o}$ of Proposition 1.6.2.
Before starting on this task, let us list some 
conventions of notation which will be used during the proof.

$\ $

{\bf 3.5 Notations.} $1^{o}$ We will use the following conventions:

-- For $1 \leq i<j \leq ks$, the complex-conjugate of the random variable
$\rho_{i,j}$ given in Proposition 1.6.2 will be denoted by $\rho_{j,i}$.
(Thus $\rho_{j,i}$ is also a random variable on $\Omega$, with values
in the unit circle.)

-- In the $2p$-tuple $\ee (1), \ldots , \ee (2p)$ which appears in the 
statement of Proposition 1.6.2, the $\ee (m)$'s which are equal to $*$ will 
be treated in algebraic expressions as if they were equal to $-1$.
(For instance ``$\sum_{b \in B} \ee (b) = 0$'', for $B$ a subset of
$\{ 1, \ldots , 2p \}$, will actually mean that 
$| \ \{ b \in B \ | \ \ee (b) = 1 \} \ |$ = 
$| \ \{ b \in B \ | \ \ee (b) = * \} \ |$. )

\vspace{10pt}

$2^{o}$ Combinatorial notations:

-- $\pp$ will denote the set of pairings of $\ppset$, where a 
pairing of $\ppset$ is as defined in Notations 1.2.1. The set of all
{\em partitions} of $\ppset$ will be denoted by $\pall$. (A partition 
$\pi = \{ B_{1}, \ldots , B_{m} \}$ of $\ppset$ is defined in the same 
way as a pairing, but without any restriction on the cardinalities of 
$B_{1}, \ldots , B_{m}$.)

-- Let $\pi$ be in $\pall$. We will say that $\pi$ is {\em r-stable}
if $r_{a} = r_{b}$ whenever $a,b \in \ppset$ belong to the same block of
$\pi$; and we will say that $\pi$ is {\em $\ee$-null}
if $\sum_{b \in B} \ee (b) = 0$ for every block $B$ of $\pi$. (Here 
$r_{a}, r_{b}$ are extracted out of the $2p$-tuple
$r_{1}, r_{2}, \ldots , r_{2p}$ appearing in the statement of 
Proposition 1.6.2, and similarly for the $\ee (b)$'s.)
Note that the index set $\ppree$ appearing in Equation (1.20) of 
Proposition 1.6.2 can be presented as
\begin{equation}
\ppree \ = \ 
\{ \ \pi \in \pp \ | \ \pi \mbox{ is $r$-stable and $\ee$-null } \}  .
\end{equation}

-- If $1 \leq i_{1}, \ldots , i_{2p} \leq ks$, then we will denote by
$ker ( i_{1}, \ldots , i_{2p} ) \in \pall$ the partition $\pi$ determined 
as follows: $a,b \in \ppset$ lie in 
the same block of $\pi$ if and only if $i_{a} = i_{b}$. 

$\ $

{\bf 3.6 Proof} {\em of part $2^{o}$ in Proposition 1.6.2.} The presentation 
of this fairly lengthy proof will be divided into several steps.

\vspace{14pt}

{\em Step 1.} The evaluation of 
$E( \ X_{r_{1}}^{\ee (1)} \cdots X_{r_{2p}}^{\ee (2p)}  \ )$
starts in the same way as the one for 
$E ( \ X_{r_{1}}^{\ee (1)} \cdots X_{r_{n}}^{\ee (n)}  \ )$
which was made in Section 3.1. We obtain the analogue of 
the Equation (3.1) of that proof:
\begin{equation}
E( \ X_{r_{1}}^{\ee (1)} \cdots X_{r_{2p}}^{\ee (2p)}  \ ) \ =
\end{equation}
\[
\frac{1}{k^{p}} \cdot \sum_{
\begin{array}{c}
{\scriptstyle 1 \leq i_{1}, \ldots , i_{2p} \leq ks \ such \ that} \\
{\scriptstyle i_{1}=r_{1}(mod \ s), \ldots , i_{2p}=r_{2p}(mod \ s) }
\end{array} }  \ \ \int_{\Omega} \varphi ( \
U_{i_{1}} ( \omega )^{\ee (1)} \cdots 
U_{i_{2p}} ( \omega )^{\ee (2p)} \ ) \ dP( \omega ).
\]

We then write the right-hand side of (3.14) as a double summation, as follows:
\[
\sum_{\pi \in \pall} \Bigl( \ \frac{1}{k^{p}} \cdot
\sum_{
\begin{array}{c}
{\scriptstyle 1 \leq i_{1}, \ldots , i_{2p} \leq ks} \\
{\scriptstyle such \ that \ \ambda (i_{1}, \ldots , i_{2p}) = \pi \ and} \\
{\scriptstyle i_{1}=r_{1}(mod \ s), \ldots , i_{2p}=r_{2p}(mod \ s) }
\end{array} }  \ \ \int_{\Omega} \varphi ( \ 
U_{i_{1}} ( \omega )^{\ee (1)} \times 
\]
\[
\times \cdots 
U_{i_{2p}} ( \omega )^{\ee (2p)} \ ) \ dP( \omega ) \ \Bigr) .
\]
In other words we write
\begin{equation}
E( \  X_{r_{1}}^{\ee (1)} \cdots X_{r_{2p}}^{\ee (2p)} \ ) \ = \ 
\sum_{\pi \in \pall} T_{\pi} ,
\end{equation}
where for $\pi \in \pall$ we set:
\begin{equation}
T_{\pi} \ := \
\frac{1}{k^{p}} \cdot
\sum_{
\begin{array}{c}
{\scriptstyle 1 \leq i_{1}, \ldots , i_{2p} \leq ks} \\
{\scriptstyle such \ that \ \ambda (i_{1}, \ldots , i_{2p}) = \pi \ and} \\
{\scriptstyle i_{1}=r_{1}(mod \ s), \ldots , i_{2p}=r_{2p}(mod \ s) }
\end{array} }  \ \ \int_{\Omega} \varphi ( \ 
U_{i_{1}} ( \omega )^{\ee (1)} \cdots 
U_{i_{2p}} ( \omega )^{\ee (2p)} \ ) \ dP( \omega ).
\end{equation}
Our strategy will be to analyze, in the following few steps of the proof,
the quantities $T_{\pi}$, $\pi \in \pall$.

\vspace{14pt}

{\em Step 2.} In this step we observe that if $\pi \in \pall$ is not 
$r$-stable (in the sense defined in Notations 3.5), then the index set of the 
summation in (3.16) is void, and hence $T_{\pi}$ = 0 (in a vacuous way).

\vspace{6pt}

{\em Proof of Step 2.}
Suppose that $\pi \in \pall$ is such that the index set of the 
summation in (3.16) is non-void. This means that there exist
$1 \leq i_{1}, \ldots , i_{2p} \leq ks$ such that
$\ambda ( i_{1}, \ldots , i_{2p} ) = \pi$ and such that
$i_{1} = r_{1} (mod \ s), \ldots , i_{2p} = r_{2p} (mod \ s)$.
Then for every $a,b$ belonging to the same block of $\pi$
we have: $i_{a} = i_{b}$ $\Rightarrow$ $r_{a}=r_{b} (mod \ s)$
(because $r_{a} = i_{a} (mod \ s)$, $r_{b}=i_{b} (mod \ s)$ )
$\Rightarrow$ $r_{a} = r_{b}$ (because $1 \leq r_{a}, r_{b}
\leq s$), and we conclude that $\pi$ is $r$-stable.

\vspace{14pt}

{\em Step 3.} Consider now a partition $\pi \in \pall$ which is $r$-stable 
but is not $\ee$-null. We show that $T_{\pi}$ = 0.

\vspace{6pt}

{\em Proof of Step 3.}
We can prove in fact a stronger statement than $T_{\pi} = 0$, namely that:
\begin{equation}
\left\{  \begin{array}{l}
\varphi ( \ 
U_{i_{1}} ( \omega )^{\ee (1)} \cdots 
U_{i_{2p}} ( \omega )^{\ee (2p)} \ ) \ = \ 0, \\
\forall \ \omega \in \Omega , \ \forall \ 1 \leq i_{1}, \ldots ,
i_{2p} \leq ks \mbox{ such that }
\ambda ( i_{1} , \ldots , i_{2p} ) = \pi .
\end{array}  \right.
\end{equation}
The proof of of (3.17) is similar to the proof of part $1^{o}$ of 
Proposition 1.6.2 (compare to Equation (3.2) in Section 3.1). Let $B$ be 
a block of $\pi$ such that $\sum_{b \in B} \ee (b) \neq 0$. If
$1 \leq i_{1}, \ldots , i_{2p} \leq ks$ are such that 
$\ambda ( i_{1}, \ldots , i_{2p} ) = \pi$, then $i_{a} = i_{b}$
for every $a,b \in B$, and it makes sense to denote by 
$i \in \{ 1 , \ldots , ks \}$
the common value of the $i_{b}$'s with $b \in B$. The commutation
relations satisfied by the unitaries $U_{1}, \ldots , U_{ks}$ give us, 
for an arbitrary $\omega \in \Omega$, an equality of the form
\begin{equation}
U_{i_{1}} ( \omega )^{\ee (1)} \cdots 
U_{i_{2p}} ( \omega )^{\ee (2p)} \ = \ c \
U_{1} ( \omega )^{\lambda_{1}} \cdots 
U_{ks} ( \omega )^{\lambda_{ks}} ,
\end{equation}
where $c$ is a constant of absolute value 1 and 
$\lambda_{1}, \ldots , \lambda_{ks} \in [-2p,2p] \cap \Z \subset 
(-L,L) \cap \Z$. The point is that 
$\lambda_{i} = \sum_{b \in B} \ee (b) \neq 0$; hence when we apply 
$\varphi$ in (3.18), we obtain 0 because of the condition (jjj-$L$) 
introduced in Section 1.6.

\vspace{14pt}

{\em Step 4.} We consider next a partition $\pi \in \pall$ which is 
$r$-stable and $\ee$-null, but is not a pairing (i.e. not all 
the blocks of $\pi$ have exactly two elements). For such a $\pi$
we prove the inequality 
\begin{equation}
| T_{\pi} | < 1/k .
\end{equation}

\vspace{6pt}

{\em Proof of Step 4.}
Observe first that the number of terms in the sum defining $T_{\pi}$ in 
(3.16) is bounded above by $k^m$, where $m$ is the number of blocks of 
$\pi$.  Indeed, constructing a $2p$-tuple $(i_{1}, \ldots , i_{2p})$
such that $\ambda ( i_{1}, \ldots , i_{2p} ) = \pi$ amounts to
constructing an injective function from the set of the blocks of
$\pi$ to the set 
$\{ 1, \ldots , ks \}$; but the requirements $i_{1}=r_{1} (mod \ s)
\ldots , i_{2p} = r_{2p} (mod \ s)$ allow only $k$ possible values
for each of the values taken by this function -- so even if the
injectivity requirement is ignored, there still are at most $k^m$
such functions which can be constructed. On the other hand, it is 
obvious that every term of the sum on the right-hand side of (3.16) 
is less or equal 1 in absolute value (contractive linear functional 
applied to a unitary).  We thus obtain that the quantity in (3.16) is
bounded in absolute value by $k^{m-p}$. But the facts that $\pi$ is
$\ee$-null and is not a pairing imply $m \leq p-1$. (Indeed, every
block of $\pi$ has an even number of elements, because $\pi$ is 
$\ee$-null; this implies $m \leq p$, with equality holding if and only
if every block of $\pi$ has exactly two elements -- which we supposed
is not the case.) Hence $k^{m-p} \leq 1/k$, and (3.19) is obtained.

\vspace{14pt}

{\em Step 5.} It is now the moment to consider a pairing 
$\pi \in \pp$, which is both $r$-stable and $\ee$-null -- or in other words,
an element $\pi \in {\cal P} ( r_{1}, \ldots , r_{2p}; \ee (1), \ldots , 
\ee (2p) )$. In this step of the proof we also fix some indices 
$1 \leq i_{1}, \ldots , i_{2p} \leq ks$ such that 
$\ambda ( i_{1}, \ldots , i_{2p} ) = \pi$ and such that 
$i_{1}=r_{1} (mod \ s), \ldots , i_{2p} = r_{2p} (mod \ s)$.
The goal of the step is to give a good approximation for the integral
\[
\int_{\Omega} \varphi ( \ 
U_{i_{1}} ( \omega )^{\ee (1)} \cdots 
U_{i_{2p}} ( \omega )^{\ee (2p)}  \ ) \ dP( \omega ) .
\]

Let us write explicitly $\pi = \{ B_{1}, \ldots , B_{p} \}$
where the blocks $B_{1}, \ldots , B_{p}$ are listed in increasing 
order of their minimal elements. The values $i_{min(B_{1})},$
$\ldots , i_{min(B_{p})} \in \{ 1, \ldots , ks \}$ are distinct,
hence it makes sense to consider the permutation
\[
\sigma \ := ord ( \ i_{min(B_{1})} , \ldots , i_{min(B_{p})}  \ )
\]
of $\pset$, which keeps track of their order ($\sigma$ defined as in
Notation 3.3). We will show that:
\[
| \ \int_{\Omega} \varphi ( \ 
U_{i_{1}} ( \omega )^{\ee (1)} \cdots 
U_{i_{2p}} ( \omega )^{\ee (2p)}  \ ) \ dP( \omega ) \ - \ 
z^{cr_{+} ( \pi , \ee , \sigma )} 
{\overline{z}}^{cr_{-} ( \pi , \ee , \sigma )}  \ |
\]
\begin{equation}
\leq  \ \frac{p(p-1)}{2} \delta .
\end{equation}

\vspace{6pt}

{\em Proof of Step 5.} Let us consider the unitaries
\[
V_{l} ( \omega ) \ := \ 
\Bigl( U_{i_{min(B_{l})}} ( \omega ) \Bigr)^{\ee ( min(B_{l}) )}, \ \
1 \leq l \leq p, \ \omega \in \Omega .
\]
Note that:
\[
\Bigl( U_{i_{max(B_{l})}} ( \omega ) \Bigr)^{\ee ( max(B_{l}) )} 
\ = \ V_{l} ( \omega )^{-1} , 
\ \ \forall \ 1 \leq l \leq p, \ \forall \  \omega \in \Omega ;
\]
this is because 
$\ee ( \mbox{max} (B_{l})) = -
\ee ( \mbox{min} (B_{l}))$
(which happens because $\pi$ is $\ee$-null), and
$i_{max(B_{l})} = i_{min(B_{l})}$ (which comes from the fact that
$\ambda ( i_{1}, \ldots , i_{2p} ) = \pi$). On the other hand for 
every $1 \leq l < m \leq p$ and every $\omega \in \Omega$ we have
the commutation relation
\[
V_{l} ( \omega ) V_{m} ( \omega ) \ = \ \Bigl( \
\rho_{i_{min(B_{l})}, i_{min(B_{m})} } ( \omega ) \ \Bigr)^{ \ee
( min(B_{l}) ) \ee ( min(B_{m}) ) } \ 
V_{m} ( \omega ) V_{l} ( \omega ),
\]
which is implied by the commutation relations known for the unitaries
$U_{i} ( \omega )$. But then the commutation 
Lemma 3.2 applies, and gives us that 
\[
U_{i_{1}} ( \omega )^{\ee (1)} \cdots
U_{i_{2p}} ( \omega )^{\ee (2p)}  
\]
\[
= \ \prod_{ \begin{array}{c}
{\scriptstyle 1 \leq l < m \leq p} \\
{\scriptstyle such \ that} \\
{\scriptstyle B_{l} \ crosses \ B_{m}}
\end{array}  } \  \Bigl( \
\rho_{i_{min(B_{l})}, i_{min(B_{m})} } ( \omega ) \ \Bigr)^{ \ee
( min(B_{l}) ) \ee ( min(B_{m}) ) } \  I.
\]

Hence we obtain:
\[
\int_{\Omega} \varphi ( \ 
U_{i_{1}} ( \omega )^{\ee (1)} \cdots 
U_{i_{2p}} ( \omega )^{\ee (2p)}  \ ) \ dP( \omega )  
\]
\[
= \ \int_{\Omega} \ \prod_{
\begin{array}{c}
{\scriptstyle 1 \leq l < m \leq p} \\
{\scriptstyle such \ that} \\
{\scriptstyle B_{l} \ crosses \ B_{m}}
\end{array}  } \  \Bigl( \
\rho_{i_{min(B_{l})}, i_{min(B_{m})} } ( \omega ) \ \Bigr)^{ \ee
( min(B_{l}) ) \ee ( min(B_{m}) ) } \ dP( \omega )
\]
\begin{equation}
= \ \prod_{ \begin{array}{c}
{\scriptstyle 1 \leq l < m \leq p} \\
{\scriptstyle such \ that} \\
{\scriptstyle B_{l} \ crosses \ B_{m}}
\end{array}  } \ \int_{\Omega} \  \Bigl( \
\rho_{i_{min(B_{l})}, i_{min(B_{m})} } ( \omega ) \ \Bigr)^{ \ee
( min(B_{l}) ) \ee ( min(B_{m}) ) } \ dP( \omega );
\end{equation}
the product and the integration could be interchanged at the last
equality sign because the random variables 
$( \rho_{ij} )_{1 \leq i<j \leq ks}$ are independent (which 
immediately implies that the random variables
$( \rho_{i_{min(B_{l})}, i_{min(B_{m})} } )_{1 \leq l<m \leq p}$
are also independent).

In the product (3.21), every factor is either within $\delta$ 
from $z$, or within $\delta$ from $\overline{z}$. In fact, one sees 
by direct inspection that:

-- if the crossing between $B_{l}$ and $B_{m}$ has positive 
orientation in $( \pi , \ee , \sigma )$, then 
\begin{equation}
| \ \int_{\Omega} \  \Bigl( \
\rho_{i_{min(B_{l})}, i_{min(B_{m})} } ( \omega ) \ \Bigr)^{ \ee
( min(B_{l}) ) \ee ( min(B_{m}) ) } \ dP( \omega ) \ - \ z \ | \ 
\leq \ \delta ;
\end{equation}

-- if the crossing between $B_{l}$ and $B_{m}$ has negative 
orientation in $( \pi , \ee , \sigma )$, then
\begin{equation}
| \ \int_{\Omega} \  \Bigl( \
\rho_{i_{min(B_{l})}, i_{min(B_{m})} } ( \omega ) \ \Bigr)^{ \ee
( min(B_{l}) ) \ee ( min(B_{m}) ) } \ dP( \omega ) \ - \ 
\overline{z} \ | \ \leq  \ \delta .
\end{equation}

In order to check (3.22-23), there are four possible cases to discuss,
according to whether $i_{min(B_{l})}$ is bigger or smaller than
$i_{min(B_{m})}$, and whether 
$\ee ( \mbox{min} (B_{l})) \cdot \ee ( \mbox{min} (B_{m}))$ is 1 or $-1$.
We show one of them, say when $i_{min(B_{l})} > i_{min(B_{m})}$ and
$\ee ( \mbox{min} (B_{l})) \cdot \ee ( \mbox{min} (B_{m})) = 1$.

The inequality $i_{min(B_{l})} > i_{min(B_{m})}$ is equivalent to 
$\sigma (l) > \sigma (m)$ (by the definition of the permutation $\sigma$
-- see (3.9) in Notation 3.3); comparing this against the formulas (1.16),
we see that $B_{l}$ and $B_{m}$ have a negative crossing. But on
the other hand:
\[
\int_{\Omega} \  \Bigl( \
\rho_{i_{min(B_{l})}, i_{min(B_{m})} } ( \omega ) \ \Bigr)^{ \ee
( min(B_{l}) ) \ee ( min(B_{m}) ) } \ dP( \omega ) 
\]
\[
= \ \Bigl( \ \int_{\Omega} \  
\overline{ \rho_{ i_{min(B_{m})}, i_{min(B_{l})} } ( \omega ) } \ 
dP( \omega ) \ \Bigr)
\]
with $i_{min(B_{m})} < i_{min(B_{l})}$; this integral is within 
$\delta$ of $\overline{z}$, by one of the hypotheses of Proposition 1.6.2.

Finally, (3.22) and (3.23) imply (3.20), via the well-known fact (easily 
checked by induction) that if $\xi_{1}, \ldots , \xi_{N}, \eta_{1}, \ldots
, \eta_{N}$ are complex numbers of value not exceeding 1, and if 
$| \xi_{1} - \eta_{1} | \leq \delta , \ldots , | \xi_{N} - \eta_{N} | 
\leq \delta$
then $| \xi_{1} \cdots \xi_{N} - \eta_{1}  \cdots \eta_{N} | \leq
N \delta$. (Here $N$ is the number of crossings of $\pi$, which cannot
exceed $p(p-1)/2$.)

\vspace{16pt}

{\em Step 6.}  In this step we fix again a pairing $\pi \in \ppree$.
We will prove the inequality:
\begin{equation}
| \ T_{\pi} \ - \ 
\frac{1}{p!} \cdot \sum_{ \sigma \in \perm} \ 
z^{ cr_{+} ( \pi , \ee , \sigma ) }
{\overline{z}}^{ cr_{-} ( \pi , \ee , \sigma ) }  \ | \ < \ 
2^{p-1} p^{2} \delta + \frac{(p+1)^{p}}{k} .
\end{equation}

\vspace{6pt}

{\em Proof of Step 6.} Let us write explicitly the partition $\pi$ fixed in 
this step as $\{ B_{1}, \ldots , B_{p} \}$, where the blocks 
$B_{1}, \ldots , B_{p}$ are listed in increasing order of their minimal 
elements. Also, let us denote:
\[
{\cal J} \ = \ \ \Bigl\{ \ (i_{1}, \ldots , i_{2p} ) \
\begin{array}{ll}
|  &  1 \leq i_{1} , \ldots , i_{2p} \leq ks,   \\
|  &  \ambda ( i_{1}, \ldots , i_{2p} ) = \pi , \\
|  &  i_{1}=r_{1} (mod \ s), \ldots , i_{2p} = r_{2p} (mod \ s)
\end{array}  \ \Bigr\} ;
\]
i.e, ${\cal J}$ is the index set of the summation defining $T_{\pi}$ in 
Equation (3.16).

For every  $2p$-tuple $(i_{1}, \ldots i_{2p}) \in {\cal J}$ we write the 
inequality (3.20) obtained in Step 5; then we sum all these inequalities.
The integrals from (3.20) will add up to $k^{p} T_{\pi}$. The terms 
``$z^{cr_{+} ( \pi , \ee , \sigma )} 
{\overline{z}}^{cr_{-} ( \pi , \ee , \sigma )}$'' from (3.20) will add up to:
\[
\sum_{ \sigma \in \perm} \ N( \sigma )
z^{ cr_{+} ( \pi , \ee , \sigma ) }
{\overline{z}}^{ cr_{-} ( \pi , \ee , \sigma ) },
\]
where for every $\sigma \in \symmp$ we denoted
\[
N( \sigma )  \ = \ | \ \{ (i_{1}, \ldots , i_{2p}) \in {\cal J} \ | \
ord( i_{ min(B_{1}) }, \ldots , i_{ min(B_{p}) } ) = \sigma \} \ | .
\]
We thus obtain, after also dividing by $k^{p}$:
\begin{equation}
| \ T_{\pi} \ - \ 
\frac{1}{k^{p}} \cdot \sum_{ \sigma \in \perm} \ N( \sigma )
z^{ cr_{+} ( \pi , \ee , \sigma ) }
{\overline{z}}^{ cr_{-} ( \pi , \ee , \sigma ) } \ | \
\leq  \ \frac{1}{k^{p}} \cdot \sum_{ \sigma \in \perm} \ N( \sigma )
\cdot \frac{p(p-1)}{2} \delta .
\end{equation}

Now, the Lemma 3.4 gives us that:
\begin{equation}
\left(  \begin{array}{c} k \\ p  \end{array} \right)
\ \leq \  N( \sigma ) \ \leq \ 
\left(  \begin{array}{c} k+p \\ p  \end{array} \right) , \ \ 
\forall \ \sigma \in \perm .
\end{equation}
One consequence of (3.26) is that the right-hand side of (3.25) is bounded 
above by:
\begin{equation}
\frac{1}{k^p} \cdot p! \cdot 
\left(  \begin{array}{c} k+p \\ p  \end{array} \right) \cdot
\frac{p(p-1)}{2} \delta \ = \ 
\frac{(k+1) \cdots (k+p)}{k^{p}} \cdot \frac{p(p-1)}{2} \delta 
\ < \ 2^{p-1} p^{2} \delta 
\end{equation}
(where at the last equality sign we used the fact that $p \leq k$).

Another consequence of (3.26) is that
\begin{equation}
| \ \frac{1}{k^{p}} \cdot \sum_{ \sigma \in \perm} \ N( \sigma )
z^{ cr_{+} ( \pi , \ee , \sigma ) }
{\overline{z}}^{ cr_{-} ( \pi , \ee , \sigma ) } 
\end{equation} 
\[
- \ \frac{1}{p!} \cdot \sum_{ \sigma \in \perm} \ 
z^{ cr_{+} ( \pi , \ee , \sigma ) }
{\overline{z}}^{ cr_{-} ( \pi , \ee , \sigma ) } \ | \ < \
\frac{(p+1)^{p}}{k} .
\]
Indeed, the left-hand side of (3.28) can be written as:
\begin{equation}
| \ \frac{1}{k^{p} p!} \cdot \sum_{ \sigma \in \perm} \ 
( p! N( \sigma ) - k^{p} ) \cdot
z^{ cr_{+} ( \pi , \ee , \sigma ) } 
{\overline{z}}^{ cr_{-} ( \pi , \ee , \sigma ) } \ | .
\end{equation} 
But for every $\sigma \in \symmp$:
\[
| \ p! N( \sigma ) - k^{p} \ | \ \leq \ \mbox{max} \Bigl( \ 
| \ p!  \left( \begin{array}{c} k \\ p \end{array} \right) - k^p \ | , \ 
| \ p!  \left( \begin{array}{c} k+p \\ p \end{array} \right) - k^p 
\ |  \ \Bigr)
\]
\[
< \ \mbox{max} \Bigl( \ k^{p} - (k-p)^{p} , \ 
(k+p)^{p} - k^{p} \ \Bigr) \ < \ k^{p-1} (p+1)^{p} ,
\]
hence the the quantity in (3.29) is dominated by
$\frac{1}{k^{p}p!} \cdot p! \cdot k^{p-1} (p+1)^{p}$ =
$\frac{(p+1)^{p}}{k}$.

The inequality (3.24) (which is the goal of Step 6) is immediately obtained
from (3.27) and (3.28).

\vspace{14pt}

{\em Step 7.} In this final part of the proof, we combine the results of 
the previous steps in order to obtain the inequality (1.20) stated in 
Proposition 1.6.2.

\vspace{6pt}

We first claim that:
\begin{equation}
| \ E( \ X_{r_{1}}^{\ee (1)} \cdots X_{r_{2p}}^{\ee (2p)}  \ )
\ - \ \sum_{ \pi \in \ppree} \ T_{\pi} \ | \ < \ \frac{(2p)!}{k} .
\end{equation}
Indeed, we know that
\[
E( \ X_{r_{1}}^{\ee (1)} \cdots X_{r_{2p}}^{\ee (2p)}  \ )
\ = \ \sum_{\pi \in \pall} T_{\pi} \ \ \mbox{(by Step 1)}
\]
\[
\ = \ \sum_{ \begin{array}{c}
{\scriptstyle \pi \in \pall ,}  \\
{\scriptstyle \pi \ r-stable} \\
{\scriptstyle and \ \ee -null}
\end{array} } T_{\pi} \ \ \mbox{(by Steps 2 and 3).}
\]
By taking into account the Equation (3.13) of Notations 3.5, we see that 
the left-hand side of (3.30) is hence equal to 
\[
| \ \sum_{ \begin{array}{c}
{\scriptstyle \pi \in \pall \setminus \pp ,}  \\
{\scriptstyle \pi \ r-stable}  \\
{\scriptstyle and \ \ee -null}
\end{array} } T_{\pi} \ | ;
\]
but by the Step 4, this is bounded above by 
$| \ \pall \setminus \pp \ | /k$, which in turn is dominated by 
$(2p)!/k$ (we used the rough estimates
$| \ \pall \setminus \pp \ | < | \ \pall \ | < (2p)!$ ). Hence (3.30)
is obtained.

We next claim that 
\[
| \ \sum_{ \pi \in \ppree } \  T_{\pi} 
\ - \ \sum_{ \pi \in \ppree } \ \frac{1}{p!} 
\cdot \sum_{ \sigma \in \perm} \ 
z^{ cr_{+} ( \pi , \ee , \sigma ) } 
{\overline{z}}^{ cr_{-} ( \pi , \ee , \sigma ) }  \ |
\]
\begin{equation}
< \ p! \cdot ( 2^{p-1} p^{2} \delta + \frac{(p+1)^{p}}{k} ).
\end{equation}
Indeed, if we write the inequality (3.24) obtained in Step 6 for every 
$\pi \in {\cal P} ( r_{1}, \ldots , r_{2p};$
$\ee (1), \ldots , \ee (2p) )$, and if we sum over 
$\pi$, we obtain that the left-hand side of (3.31) is bounded above by
\[
| \ \ppree \ | \cdot ( 2^{p-1} p^{2} \delta + \frac{(p+1)^{p}}{k} ).
\]
The latter quantity is in turn dominated by 
$p! \cdot ( 2^{p-1} p^{2} \delta + \frac{(p+1)^{p}}{k} )$, because the
number of pairings in $\pp$ which are $\ee$-null (but not necessarily
$r$-stable) is exactly $p!$ .

The desired inequality (1.20) immediately follows from (3.30), (3.31),
and the rough estimates $(2p)!+ p!(p+1)^{p} < (2p+1)!$,
$p! 2^{p-1} p^{2} < (2p+1)!$. {\bf QED}

$\ $

For the rest of the section we move to the framework of Proposition
1.6.3. The proof of 1.6.3 is in many respects similar to the one of 1.6.2.
For this reason we will not write the arguments in the same detail, and
occasionally we will leave it as an exercise to the reader to check that 
parts of the proof of 1.6.2 can be trivially adjusted to the current 
situation.

$\ $

{\bf 3.7 Proof} {\em of Proposition 1.6.3.} Let $E: \BB \rightarrow \C$ be
the linear functional defined as in Equation (1.12) of Notation 3.4. We fix 
$n \geq 1$, and $\rs$, $\ees$, about which we will show that the limit
\begin{equation}
\lim_{k \rightarrow \infty} \ E( \ 
X_{r_{1};k}^{\ee (1)} \cdots X_{r_{n};k}^{\ee (n)} \ )
\end{equation}
exists and is equal to the right-hand side of Equation (1.17).

In connection to these $n$, $r_{1}, \ldots , r_{n}$, 
$\ee (1), \ldots , \ee (n)$ that are fixed, we will use combinatorial 
notations similar to some of those set in Notations $3.5.2^{o}$:
$\palln$ will denote the set of 
all the partitions of $\{ 1, \ldots , n \}$; and we will say that
$\pi \in \palln$ is ``$r$-stable'' if $r_{a} = r_{b}$ whenever 
$a,b \in \{ 1, \ldots , n \}$ belong to the same block of $\pi$.

We leave it as an exercise to the reader to verify that the Step 1 of the 
proof in Section 3.6 can be performed in the current situation, and leads
to the following analogue of the Equations (3.15-16):
\begin{equation}
E( \ X_{r_{1};k}^{\ee (1)} \cdots X_{r_{n};k}^{\ee (n)} \ ) \ = \ 
\sum_{\pi \in \palln} T_{\pi , k}, \ \ \forall \ k \geq 1,
\end{equation}
where for $\pi \in \palln$ we set:
\begin{equation}
T_{\pi , k} \ := \
\frac{1}{k^{n/2}} \cdot
\sum_{
\begin{array}{c}
{\scriptstyle 1 \leq i_{1}, \ldots , i_{n} \leq ks} \\
{\scriptstyle such \ that \ \ambda (i_{1}, \ldots , i_{n}) = \pi \ and} \\
{\scriptstyle i_{1}=r_{1}(mod \ s), \ldots , i_{n}=r_{n}(mod \ s) }
\end{array} }  \ \ \int_{\Omega} \varphi ( \ 
U_{i_{1}} ( \omega )^{\ee (1)} \cdots 
U_{i_{n}} ( \omega )^{\ee (n)} \ ) \ dP( \omega ).
\end{equation}

It is also clear that the Step 2 of the proof in Section 3.6 can be repeated
identically, and leads to the conclusion that $T_{\pi , k} = 0$ for every 
$k \geq 1$ and every $\pi \in \palln$ which is not $r$-stable. Thus the 
partitions which are not $r$-stable can be ignored in the summation on the 
right-hand side of (3.33).

We next observe that for every $k \geq 1$ and every $\pi \in \palln$ which 
is not $r$-stable, we have the inequality:
\begin{equation}
| T_{\pi ,k} | \ \leq \ k^{-(\frac{n}{2} - | \pi |)} ,
\end{equation}
where $| \pi |$ stands for the number of blocks of the partition $\pi$.
The verification of this inequality is very similar to the argument shown 
in Step 4 of the proof in 3.6, and is left to the reader. Due to the
fact that in (3.33) we are actually interested only in what happens when 
$k \rightarrow \infty$, the inequality (3.35) shows that in the summation
on the right-hand side of (3.33) we can also safely ignore all the 
$r$-stable partitions $\pi$ such that $| \pi | < n/2$.

Now let us remark that $T_{\pi , k} = 0$ for every $k \geq 1$
and for every $\pi \in \palln$ which has at least one singleton (i.e. a 
block with one element). Indeed, let us suppose that the partition $\pi$
has a one-element block $B = \{ b \}$, $1 \leq b \leq n$. Then for every 
$1 \leq i_{1}, \ldots , i_{n} \leq ks$ such that 
$\ambda ( i_{1}, \ldots , i_{n} ) = \pi$, the monomial 
$U_{i_{1}} ( \omega )^{\ee (1)} \cdots U_{i_{n}} ( \omega )^{\ee (n)}$ is 
brought by the commutation relations (1.21) to the form
$c \ U_{1} ( \omega )^{\lambda_{1}} \cdots U_{ks} ( \omega )^{\lambda_{ks}}$,
where $|c| =1$, $\lambda_{1}, \ldots , \lambda_{ks} \in \Z$, and -- most 
importantly here -- $\lambda_{i_{b}} = \pm 1$. But then the hypothesis (d)
of Proposition 1.6.3 gives that $\varphi 
(U_{i_{1}} ( \omega )^{\ee (1)} \cdots U_{i_{n}} ( \omega )^{\ee (n)} )$ =
0, and the equality $T_{\pi , k} = 0$ follows.

The conclusion of the preceding three paragraphs is that in the summation on 
the right-hand side of (3.33) we may keep (without affecting what happens 
when $k \rightarrow \infty$) only the terms which correspond to partitions 
$\pi \in \palln$ that are $r$-stable, satisfy $| \pi | \geq n/2$, and have
no singletons.

However, if $n$ is odd, then there are no partitions at all which satisfy 
$| \pi | \geq n/2$ and at the same time have no singletons.
This simply means that if $n$ is odd, then the limit in (3.32) exists and 
is equal to 0 (and the case of odd $n$ is thus settled).

If $n$ is even, it is immediate that a partition $\pi \in \palln$ satisfies
$| \pi | \geq n/2$ and has no singletons if and only if it is a pairing.
Thus in the case of even $n$, the summation on the right-hand side of (3.33)
can be restricted to the set of $r$-stable pairings of
$\{ 1, \ldots , n \}$.

From now on and until the end of the proof we will assume that $n$ is even,
$n = 2p$ with $p$ positive integer. Similarly to the terminology introduced
in the Notations $3.5.2^{o}$, we will say that a pairing $\pi$ =
$\{ \ \{ a_{1}, b_{1} \}, \ldots , \{ a_{p}, b_{p} \} \ \}$ of 
$\{ 1, \ldots , n \}$ is $\ee$-null if $\ee ( a_{i} ) \neq \ee ( b_{i} )$,
$\forall \ 1 \leq i \leq p$ (where $\ees$ are as fixed at the beginning of
the proof). By taking into account the conclusion of the preceding
paragraph, and by examining at the same time the right-hand side of Equation
(1.17), we see that the proof will be completed if we can show that:
\begin{equation}
\lim_{k \rightarrow \infty} T_{\pi , k} \ = \ 0
\end{equation}
for every pairing $\pi$ of $\{ 1, \ldots , 2p \}$ which is $r$-stable but 
not $\ee$-null; and
\begin{equation}
\lim_{k \rightarrow \infty} T_{\pi , k} \ = \ 
\frac{1}{p!} \sum_{\sigma \in \symmp} \
\sum_{\pi \in \ppree} \ z^{cr_{+} ( \orpi )} \cdot 
{\overline{z}}^{cr_{-} ( \orpi )} 
\end{equation}
for every pairing $\pi$ of $\{ 1, \ldots , 2p \}$ which is both $r$-stable 
and $\ee$-null.

The limit in (3.36) holds trivially: $T_{\pi , k} = 0$ for every $k \geq 1$
and every pairing $\pi$ which is $r$-stable but not $\ee$-null. This is a
direct application of the hypothesis (d) in Proposition 1.6.3, and is 
left to the reader. (The discussion is similar to the one which ruled out 
the partitions with singletons, but this time one uses the case when there
exists a $\lambda_{j}$ equal to $\pm 2$.)

So it suffices if from now on we fix a pairing $\pi$ of 
$\{ 1, \ldots , 2p \}$ which is both $r$-stable and $\ee$-null, and we prove
that the limit (3.37) holds. We denote the quantity on the right-hand side 
of (3.37) by $Q_{\pi}$. We will also fix a number $\beta > 0$, and we will
show that $| T_{\pi , k} - Q_{\pi} | < \beta$ if $k$ is sufficiently large.

Denote $\delta := \beta / ( 2^{p}p^{2} )$ and $L := n+1$. By the hypotheses
(b) and (c) of Proposition 1.6.3, there exists $m_{o} \geq 1$ such that 
for every $m_{o} \leq m<n$ we have that
$| z - \int_{\Omega} \rho_{m,n} | \leq \delta$, and that 
$U_{m},U_{m+1}, \ldots , U_{n}$ is an $L$-mimic of a Haar family. We fix 
$k_{o}$ such that $k_{o}s+1 \geq m_{o}$. For every $k > k_{o}$ we will write
$T_{\pi ,k}$ as a sum,
\begin{equation}
T_{\pi , k} \ = \ T_{\pi , k} '  + T_{\pi , k} '',
\end{equation}
by splitting the index set of the sum in (3.34), which defines 
$T_{\pi , k}$, into two disjoint parts: in $T_{\pi , k} '$ we take the terms
indexed by $n$-tuples $(i_{1}, \ldots , i_{n})$ such that 
$k_{o}s+1 \leq i_{1}, \ldots , i_{n} \leq ks$, and in $T_{\pi , k} ''$ 
we take the rest of the terms (indexed by $n$-tuples 
$(i_{1}, \ldots , i_{n})$ such that 
$\min (i_{1}, \ldots , i_{n}) \leq k_{o}s$).

Note that for every $k>k_{o}$, the random variables 
$( \rho_{i,j} )_{k_{o}s+1 \leq i<j \leq ks}$ and the random unitaries 
$( U_{i} )_{i=k_{o}s+1}^{ks}$ fall under the hypotheses of Proposition 1.6.2
(for the chosen values of $\delta$ and $L$); thus the estimates found in the 
proof of Proposition 1.6.2 apply to this situation. Out of these estimates,
the one which we need here is the inequality (3.24) established in the 
Step 6 of Section 3.6. When reporting to the current notations,
``$T_{\pi}$'' of (3.24) has to be replaced by:
\[
\frac{1}{(k-k_{o})^{p}} \cdot
\sum_{
\begin{array}{c}
{\scriptstyle k_{o}s+1 \leq i_{1}, \ldots , i_{2p} \leq ks} \\
{\scriptstyle such \ that \ \ambda (i_{1}, \ldots , i_{2p}) = \pi \ and} \\
{\scriptstyle i_{1}=r_{1}(mod \ s), \ldots , i_{2p}=r_{2p}(mod \ s) }
\end{array} }  \ \ \int_{\Omega} \varphi ( \ 
U_{i_{1}} ( \omega )^{\ee (1)} \cdots 
U_{i_{2p}} ( \omega )^{\ee (2p)} \ ) \ dP( \omega );
\]
but this is exactly $k^{p}/(k-k_{o})^{p} T_{\pi , k} '$, with 
$T_{\pi , k} '$ taken from (3.38). So the inequality (3.24) becomes in
this situation:
\[
\Bigl| \ \left( \frac{k}{k-k_{o}} \right)^{p} T_{\pi , k}' - Q_{\pi} \
\Bigr| \ < \ 2^{p-1}p^{2} \delta + \frac{ (p+1)^{p} }{k-k_{o}} ;
\]
or after multiplication with $(k-k_{o})^{p}/k^{p} < 1$, and after taking 
into account the relation between $\beta$ and $\delta$:
\begin{equation}
\Bigl| \ T_{\pi , k}' - \left( \frac{k-k_{o}}{k} \right)^{p} Q_{\pi} \
\Bigr| \ < \ \frac{\beta}{2} + \frac{ (p+1)^{p} }{k-k_{o}} , \ \ 
\forall \ k> k_{o}.
\end{equation}

But on the other hand, a counting argument very similar to the one shown 
in Step 4 of Section 3.6 shows that for $k > k_{o}+p$ there are less than
$k^{p}$ terms in the summation defining $T_{\pi ,k}$, and there are more
than $(k-k_{o}-p)^{p}$ terms in the summation defining $T_{\pi ,k} '$; this
implies that there are less than $k^{p} - (k-k_{o}-p)^{p}$ terms in the
summation defining $T_{\pi ,k} ''$, and consequently that:
\begin{equation}
| \ T_{\pi , k}'' | \ < \ \frac{ k^{p}-(k-k_{o}-p)^{p} }{ k^{p} } , \ \ 
\forall \ k> k_{o}+p.
\end{equation}

Finally, for $k>k_{o}+p$ we can write:
\[
| \ T_{\pi , k} - Q_{\pi} \ | \ \leq \
\Bigl| \ T_{\pi , k}' - \left( \frac{k-k_{o}}{k} \right)^{p} 
Q_{\pi} \ \Bigr| + 
\Bigl| \ 1 - \left( \frac{k-k_{o}}{k} \right)^{p} \ \Bigr| 
\cdot \Bigl| \ Q_{\pi} \ \Bigr| + 
\Bigl| \ T_{\pi , k} '' \ \Bigr|
\]
\[
< \ \frac{\beta}{2} + \frac{ (p+1)^{p} }{ k-k_{o} } +
\Bigl| \ 1 - \left( \frac{k-k_{o}}{k} \right)^{p} \ \Bigr| \cdot 
\Bigl| \ Q_{\pi} \ \Bigr| + 1 - 
\left( \frac{ k-k_{o}-p }{ k } \right)^{p} 
\]
( by Equations (3.39) and (3.40) ), and the latter expression is clearly
smaller than $\beta$ if $k$ is large enough. {\bf QED}

$\ $

$\ $

{\bf Acknowledgement:} The second-named author acknowledges the hospitality
of the Henri Poincar\'e Institute (Centre Emile Borel -- UMS 839 IHP
CNRS/UPMC) in Paris, where he visited during the final stage of the 
preparation of this paper.

$\ $

$\ $

\end{document}